\input amstex
\magnification 1200
\documentstyle{amsppt}
\NoRunningHeads
\NoBlackBoxes

\topmatter
\title 
An analogue of the field-of-norms functor  and the Grothendieck Conjecture
\endtitle

\author Victor Abrashkin 
\endauthor

 

\abstract The paper contains a construction of an analogue of 
the Fontaine-Wintenberger field-of-norms functor for 
higher dimensional local fields. This construction 
is done completely in terms of the ramification 
theory of such fields. It is applied to deduce the mixed characteristic case of a local analogue 
of the Grothendieck Conjecture for these fields from its characteristic 
$p$ case, which was proved earlier by the author. 
\endabstract
\endtopmatter

\def\pr{\operatorname{pr}}
\def\Ker{\operatorname{Ker}}
\def\ab{\operatorname{ab}}
\def\Gal{\operatorname{Gal}}
\def\sep{\operatorname{sep}}
\def\Hom{\operatorname{Hom}}
\def\id{\operatorname{id}}

\def\Aut{\operatorname{Aut}}

\def\LF{\operatorname{LF}}
\def\LC{\operatorname{LC}}
\def\char{\operatorname{char}}
\def\m{\operatorname{m}}
\def\cl{\operatorname{cl}}
\def\RLF{\operatorname{RLF}}
\def\SAPF{\operatorname{SAPF}}
\def\char{\operatorname{char}}
\def\Fil{\operatorname{Fil}}

\document 

\subhead 0. Introduction 
\endsubhead 
\medskip 

The field-of-norms functor [FW1,2] allows to identify 
the Galois groups of some infinite extensions 
of $\Bbb Q_p$ with those of complete discrete valuation fields of characteristic $p$. 
This functor is an essential component of Fontaine's theory of 
$\varphi $-$\Gamma $-modules --- one of most powerful tools in the modern study 
of $p$-adic representations cf. e.g.[Ber]. Other areas of very impressive applications are the Galois 
cohomology of local fields [He], 
arithmetic aspects of dynamical systems [LMS], 
explicit reciprocity formulae [Ab2,3], [Ben], 
a description of the structure of ramification filtration [Ab7], 
the proof of an analogue of the Grothendieck Conjecture for 1-dimensional local fields [Ab4].

A local analogue of the Grothendieck Conjecture establishes an opportunity to recover the structure 
of a local field from the structure of its absolute Galois group together with its 
filtration by ramification subgroups. The study of this conjecture in the context 
of higher dimensional local fields became actual due to a recent development of 
ramification theory for such fields [Zh2], [Ab5]. 
The characteristic $p$ case of the Conjecture has been already considered in [Ab6]. 
(Notice that the restriction to 2-dimensional fields is not essential in [Ab6] ---  
the method works for any dimension $N\geqslant 2$.) 
This result could lead to the proof of the mixed characteristic case of that conjecture 
if there were a suitable analogue of the field-of-norms functor for 
higher dimensional local fields. 

The construction of such a functor is suggested in the present paper. 
In our setting we replace the appropriate category of infinite extensions of $\Bbb Q_p$ 
by the category $\Cal B^a(N)$ of infinite increasing field towers 
$K_0\subset K_1\subset \dots \subset K_n\subset\dots $\ \  with restrictions on the upper ramification numbers 
of the intermediate extensions $K_{n+1}/K_{n}$. In order to introduce the set of elements of 
the corresponding field-of-norms one can't use in such towers the sequences of norm compatible 
elements   
but it is still possible to work with the sequences of elements $a_n\in \Cal O_{K_n}$ such that 
$a_n\equiv a_{n+1}^p\operatorname{mod}p^{c}$, where $0<c\leqslant 1$ is independent on $n$. 

The main difficulty in the realization of this idea comes from the fact that the construction 
of ramification theory for an $N$-dimensional local 
field $L$ depends on the choice of its $F$-structure, i.e. on the choice of 
the subfields $L(i)$ of $i$-dimensional constants, 
where $1\leqslant i\leqslant N$. On the other hand, in order to work with elements of $L$ 
one should use one or another choice of its local parameters. This choice can be made compatible with a 
given $F$-structure only after passing to some 
finite \lq\lq semistable\rq\rq\ extension of $L$. This explains why we have a precise analogue of 
the Fontaine-Wintenberger functor only for a subcategory of \lq\lq special\rq\rq\ towers 
$\Cal B^{fa}(N)$ in $\Cal B^a(N)$. Nevertheless, the construction of the functor 
can be extended to the whole category $\Cal B^a(N)$ and can be applied to deduce 
the mixed characteristic case of the Grothendieck Conjecture 
from the characteristic $p$ case. 

We now briefly explain the content of the paper. 

Section 1 contains preliminaries: definitions and simplest properties of 
\linebreak 
$N$-dimensional local fields $L$. We pay the special attention to the concept 
of $P$-topology --- this is the topology on $L$, which  accumulates properties of 
$N$ valuation topologies which can be attached to $L$. Then the Witt-Artin-Schreier duality 
and the Kummer theory allow us to transfer the $P$-topological structure to the group 
$\Gamma _L^{\ab }(p)$, where $\Gamma _L(p)$ is the Galois group of the 
maximal $p$-extension of $L$. This structure gives an opportunity  
to work  with $\Gamma _L(p)$ in terms of generators, cf. [Ab6].  

Section 2 contains a \lq\lq co-analogue\rq\rq\  of Epp's elimination wild ramification. 
This statement deals with a subfield of $(N-1)$-dimensional 
constants in an $N$-dimensional local field. (Most widely known 
interpretation of Epp's procedure 
deals with a subfield of 1-dimensional constants.) 
Our proof establishes an elimination procedure which is similar to the procedure  
developed in [ZhK], where it was shown that an essential part of 
such elimination can 
be done inside a given deeply ramified extension in the sense of [CG].  
This elimination procedure is required to justify the main starting point in the construction of the 
ramification theory for higher dimensional local fields from [Ab5]. 
(The original arguments from [Ab5] were not complete, cf. remark in n.2.1.)

Section 3 contains a brief introduction into the ramification theory and 
contains a version of Krasner's Lemma in the context of higher 
dimensional local fields. In Section 4 we introduce and study 
the categories of special towers $\Cal B^a(N)$ and $\Cal B^{fa}(N)$. 
These towers play a role of strict arithmetic profinite extensions from  
the Fontaine-Wintenberger construction of the field-of-norms functor. 

In section 5 we explain the construction of a family of characteristic $p$ 
local fields $X(K_{\centerdot })$, where $K_{\centerdot }\in\Cal B^{fa}(N)$,  
and prove that all such fields can be identified after (roughly speaking) 
taking inseparable extensions of constant subfields of lower dimension. 
These fields will play a role of the field-of-norms attached to a given 
tower $K_{\centerdot }\in\Cal B^{fa}(N)$. In section 6 we apply 
Krasner's Lemma from section 3 to establish all expected properties 
of the correspondence $K_{\centerdot }\mapsto \Cal K\in X(K_{\centerdot })$, 
where $K_{\centerdot }\in\Cal B^{fa}(N)$. 
In section 7 we use these properties to define the analogue 
$\Cal X_{K_{\centerdot }}$, where $K_{\centerdot }\in\Cal B^{fa}(N)$, of 
the field-of-norms functor. 
In addition, we use the operation of 
radical closure to extend this construction  
to the whole category $\Cal  B^a(N)$. In section 8 we prove that the identification 
of the Galois groups $\Gamma _{\widetilde{K}}$ (where $\widetilde{K}$ is the $p$-adic 
closure of the composite of all fields from the tower $K_{\centerdot }$)  
and $\Gamma _{\Cal K}$ becomes $P$-continuous when being 
restricted to their maximal abelian $p$-quotients. The proof is based 
on a higher dimensional version of the relation 
between the Witt-Artin-Schreier theory for $\Cal K$ and the Kummer theory for 
$\widetilde{K}$ from [Ab2]. This relation  
together with the proof of the compatibility of the proposed 
field-of-norms functor with the class field theories for 
$\Cal K$ and $\widetilde{K}$, leads to another proof 
of the explicit reciprocity formula from [Vo] (cf. also [Ka]) 
---  
the details will appear later elsewhere. 

Finally, the $P$-continuity result from n.8 allows us   
to prove in section 9 the mixed characteristic case 
of the Grothendieck Conjecture. Notice that the construction 
of the higher dimensional version of the field-of-norms functor 
from this paper 
is especially adjusted to the proof of this conjecture and was motivated 
by Deligne's paper [De]. It should be 
also mentioned that there are definite ideological links with 
methods of the paper [Fu], where 
the construction of Coleman power series was developed in the context of 
2-dimensional local fields with further applications to the construction 
of $p$-adic $L$-functions.

The author is very grateful to the MPI (Bonn) for hospitality, where 
a part of this paper was being written. 

\medskip 
\medskip

 \subhead 1. Preliminaries 
\endsubhead 
\medskip 

\subsubhead 1.1. The concept of higher dimensional local field 
\endsubsubhead 

Let $K$ be an $N$-dimensional local field, where $N\in\Bbb Z_{\geqslant 0}$. 
In other words, if $N=0$ then $K$ is a finite field and for $N\geqslant 1$, 
$K$ is a complete discreet valuation field with residue field $K^{(1)}$ 
which is an $(N-1)$-dimensional local field. We use the notation $K^{(N)}$ 
for the last residue field of $K$. 

Let $O_K^{(1)}$ be the valuation ring of $K$ with respect to first valuation 
and let $\alpha :O_K^{(1)}\longrightarrow K^{(1)}$ be a natural projection. 
Define the valuation ring $\Cal O_K$ of $K$ by setting for $N=0$, 
$\Cal O_K=K$ and for $N\geqslant 1$, $\Cal O_K=\alpha ^{-1}(\Cal O_{K^{(1)}})$. 
Recall that a system $t_1,\dots ,t_N\in \Cal O_K$ is a system of local 
parameters in $K$ if $t_1$ is a uniformiser in $O_K^{(1)}$ and 
$\alpha (t_2),\dots ,\alpha (t_N)$ is a system of local parameters in $K^{(1)}$.

In terms of such system of local parameters any element $\xi \in K$ can be uniquely 
presented as a power series of the following form

$$\xi =\sum\Sb \bar a=(a_1,\dots ,a_N) \endSb \alpha _{\bar a}t_1^{a_1}\dots t_N^{a_N}.$$ 
Here all coefficients $\alpha _{\bar a}$ are either elements of $K^{(N)}$ if $\char K=p>0$ 
or the Teichmuller representatives of those if $\char K=0$. All indices $a_i\in\Bbb Z$ and there are 
integers (which depend on $\xi $) $A_1$, $A_2(a_1)$, 
\dots , 
$A_N(a_1,\dots ,a_{N-1})$ such that  
$\alpha _{\bar a}=0$ if either 
$a_1<A_1$, or $a_2<A_2(a_1)$,\dots , or $a_N<A_N(a_1,\dots ,a_{N-1})$. 

There is an important concept of $P$-topology on $K$ which brings 
into correlation all $N$ valuation topologies related to $K$. 
The $P$-topological structure provides us with 
a reasonable treatment of morphisms of higher dimensional 
local fields. 
We discuss this structure briefly in n.1.2 below. Notice that if $f:K\longrightarrow L$ is a 
$P$-continuous 
 morphism of higher dimensional local fields 
then $E=f(K)$ is a closed subfield 
in $L$ (i.e. $O^{(1)}_E$ is closed in $O^{(1)}_L$ with respect to 
first valuation and $E^{(1)}$ is closed in $L^{(1)}$), for any system $t_1,\dots ,t_N$ 
of local parameters in $K$ their images $f(t_1),\dots ,f(t_N)$ are local parameters in $E$ 
and their knowledge determines the morphism $f$ uniquely.

Our considerations will be limited with local fields $K$ such that 
$\char K^{(1)}=p$ where $p$ is a fixed prime number (such fields possess  
most interesting arithmetic structure).  
Under this assumption there is the following classification of $N$-dimensional local fields: 

--- if $\char K=p$, then $K=k((t_N))\dots ((t_1))$ where $k=K^{(N)}$ is the last residue field of $K$. 
As a matter of fact, this result is equivalent to the 
existence of a system of local parameters 
$t_1,\dots ,t_N$ in $K$. 

--- if $\char K=0$, then $K\supset \Bbb Q_p$ and we can introduce a canonical subfield 
$K(1)$ of 1-dimensional constants in $K$: this is the algebraic closure of 
$\Bbb Q_p$ in $K$. Suppose a uniformising element $t_1$ of $K(1)$ can be 
included in a system of local parameters $t_1,t_2,\dots ,t_N$ of $K$. 
Then $K=K(1)\{\{t_{N}\}\}\dots \{\{t_2\}\}$ and such $K$ is called standard. 
 Otherwise, there is a finite extension 
$E$ of $K(1)$ such that the composite $KE$ is standard. 

The above result about the characteristic 0 fields 
is implied by the following version of Epp's theorem [Epp], which holds for all 
(not necessarily characteristic 0) 
higher dimensional local fields $K$: 

--- {\it suppose $K$ is an $N$-dimensional field and $K(1)$ is its 
subfield of 1-dimensional constants; then there is a finite extension $E$ of $K(1)$ such 
that the fields $KE$ and $E$ have a common uniformising element (with respect to the first 
valuation in $K$).} 
\medskip 

\subsubhead {\rm 1.2.} The concept of $P$-topology 
\endsubsubhead 

Let $K$ be an $N$-dimensional local field. Its $P$-topology can be described explicitly 
by induction on $N$ in terms of any chosen system $t_1,\dots ,t_N$ of local 
parameters of $K$ by constructing a basis of open 0-neighborhoods $\Cal U_b(K)$, cf. [Zh1].  
We shall consider the following three cases:
\newline 
a) $\char K=p$;
\newline 
b) $\char K=0$, $\char K^{(1)}=p$ and $t_1$ is a local parameter in  
$K(1)$;
\newline 
c) $K$ is a finite extension of $E$, which satisfies the above assumptions from b).
\medskip 

{\it The case } a).

Here $K=k((t_N))\dots ((t_1))$, where $k$ is a finite field of 
characteristic $p$. If $N=0$ then $\Cal U_b(K)$ contains by definition only 
one set $\{0\}$. Then the family $\Cal U(K)$ of all open sets in $K$ consists of 
all subsets of $K$. 
Suppose $N\geqslant 1$. Let $\bar t_N$,\dots ,$\bar t_2$ be the images of 
$t_N$,\dots , $t_2$ in $K^{(1)}$. Then $K^{(1)}=k((\bar t_N))\dots ((\bar t_2))$ 
and we can use the correspondences $\bar t_N\mapsto t_N$, \dots ,$\bar t_2\mapsto t_2$ 
and $\alpha\mapsto\alpha $ for $\alpha\in k$, 
to define the embedding $h:K^{(1)}\longrightarrow K$. Then 
$\Cal U_b(K)$ consists of the sets $\sum\Sb a\in\Bbb Z \endSb t_1^ah(U_a)$, 
where all $U_a\in\Cal U(K^{(1)})$ and $U_a=K^{(1)}$ for $a\gg 0$. 
\medskip 

{\it The case } b). 

Here again the images $\bar t_2$,\dots , $\bar t_N$ give a system of local parameters of 
$K^{(1)}$ and the family of all open subsets of $K^{(1)}$ is already defined by induction. 
So, we use again the map $h:K^{(1)}\longrightarrow K$, which is 
determined by the correspondences $\bar t_i\mapsto t_i$, $i=2,\dots ,N$, and $\alpha \mapsto [\alpha ]$ 
for $\alpha\in k$, 
and proceed along lines in the case a).
\medskip 

{\it The case } c). 

If $[K:E]=n$, then the $P$-topological structure on $K$ 
comes from any isomorphism of $E$-vector spaces $K\simeq E^n$ 
and the $P$-topological structure on $E$. 
\medskip 

It is well-known that $K$ is an additive $P$-topological group but 
the multiplication in $K$ has very bad $P$-topological 
properties. Later we need to study the $P$-continuity of 
maps between objects obtained from $K$-spaces by duality. 
For this reason we shall use the following description 
of compact subsets in $K$. 

Introduce a basis $\Cal C_b(K)$ of compact subsets in $K$. 
In other words, if $\Cal C_b(K)$ is such a family then 
any compact subset $D$ in $K$ will appear as 
a closed subset of some $C\in \Cal C_b(K)$. Proceed again 
by induction on the dimension $N$ of $K$ according to above 
assumptions a)-c) about $K$. 

In the case a) $\Cal C_b(K)$ will consist of the only one set $\{K\}$ if $N=0$. 
If $N\geqslant 1$ then in the cases a) and b) we can    
use the map $h:K^{(1)}\longrightarrow K$ to define 
$\Cal C_b(K)$ as the family of subsets $\sum\Sb a\in\Bbb Z \endSb t_1^ah(C_a)$, 
where $C_a\in\Cal C_b(K)$ and $C_a=\{0\}$ for $a\ll 0$. In the case c), we just set 
$\Cal C_b(K)=\{C^n\ |\ C\in\Cal C_b(E)\}$. 

\proclaim{Proposition 1.1} The above defined family $\Cal C_b(K)$ is 
a basis of $P$-compact subsets in $K$.
\endproclaim 

\demo{Proof} Proceed by induction on $N$ when 
$K$ satisfies the assumptions from the cases a) and b). 
The case $N=0$ is clear. 

Let $N\geqslant 1$. 
Prove first that $\Cal C_b(K)$ consists of compact subsets in $K$. 
Suppose $C=\sum t_1^ah(C_a)\in\Cal C_b(K)$. Notice 
first, that each $h(C_a)$ is $P$-compact in $K$. For any $b\in\Bbb Z$, 
set $C_{\leqslant b}=\sum\Sb a\leqslant b\endSb t_1^ah(C_a)$. Then 
$C_{\leqslant b}$ is $P$-homeomorphic to the product of 
finitely many compact sets $h(C_a)$, $a\leqslant b$. Therefore, 
$C_{\leqslant b}$ is $P$-compact. Finally, 
$C=\mathbin{\underset{b}\to\varprojlim}C_{\leqslant b}$ as 
$P$-topological sets. So, $C$ is compact. 

Suppose $D$ is a $P$-compact subset in $K$. Take $a_0\in\Bbb Z$ such that 
\linebreak 
$D\subset \sum\Sb a\geqslant a_0\endSb t_1^ah(K^{(1)})$  
(such $a_0$ exists because $D$ is compact).  
From the deinition of the $P$-topology it follows that all projections 
$\pr _a:D\longrightarrow K^{(1)}$ 
(where for any $d\in D$, $d=\sum t_1^ah(\pr _a(d))$) 
are open maps. Therefore, all $\pr _a(D)$ are compact subsets in 
$K^{(1)}$. By induction there are $C_a\in\Cal C_b(K^{(1)})$ 
such that $\pr _a(D)$ are closed subsets in $C_a$. So, 
$D$ is a subset in the $P$-compact set $\sum t_1^ah(C_a)\in\Cal C_b(K)$. 

Finally, the case c) follows from the definition of 
the $P$-topology as the product topology associated with the 
$P$-topology on $E$.
The proposition is proved.
\enddemo 

The following proposition can be proved easily by induction on $N$. 

\proclaim{Proposition 1.2} For any $C_1,C_2\in\Cal C_b(N)$, 
$C_1+C_2\in\Cal C_b(K)$ and $C_1C_2\in\Cal C_b(K)$.
\endproclaim 

\remark{Remark} A small modification of the above arguments proves the existence 
of a base of compact subsets $\Cal C'_b(K)$, which consists of additive subgroups of $K$.
\endremark 
\medskip

 \subhead 2. Higher dimensional elimination of wild ramification 
\endsubhead 
\medskip

2.1.    
 Introduce the category $\LC $ of higher dimensional local fields with a 
given subfield of constants of codimension 1. The objects in $\LC $ are 
couples $(K,E)$ where $K$ is a local field of dimension $N\geqslant 1$ 
and $E$ is a topologically closed subfield of dimension $N-1$ which is algebraically closed in $K$. 
If $N=1$ and $\char K=0$ we shall agree by definition to 
take as $E$ the maximal unramified extension of $\Bbb Q_p$ in $K$, i.e. in this case a  
1-dimensional field will play a role of a subfield of 
0-dimensional constants.  
Morphisms $(K,E)\longrightarrow (K',E')$ in the category $\LC $ 
are given by $P$-continuous morphisms 
of local fields $f:K\longrightarrow K'$ such that $f(E)\subset E'$. 

We shall use the notation $\LC (N)$ for the full subcategory in $\LC $ 
consisting of $(K,E)$, where $K$ is an $N$-dimensional field. Notice that 
$\LC (1)$ is equivalent to the usual category of complete discrete valuation fields 
with finite residue field of characteristic $p$. 

\remark{Remark} Suppose $(K,E)\in\LC $. Then 
there is a natural embedding of first residue fields 
$E^{(1)}\subset K^{(1)}$ but $(K^{(1)},E^{(1)})$ 
is not generally an object of the category 
$\LC (N-1)$,  because $E^{(1)}$ is not 
generally algebraically closed in $K^{(1)}$. 
Notice that it is separably closed in $K^{(1)}$:  
otherwise, $E$ will possess a non-trivial unramified extension in $K$. 
\endremark 
\medskip 

\definition{Definition} $(K,E)\in\LC (N)$ is standard if there is a system of local parameters 
$t_1,\dots ,t_N$ in $K$ such that $t_1,\dots ,t_{N-1}$ is a system of local 
parameters in $E$. 
In other words, 
if $(K,E)$ is standard then 
there is a $t_N\in K$ which extends 
any system of local parameters in $E$ to a system of local parameters in $K$. 
Such an element $t_N$ of $K$ will be called an $N$-th local parameter 
in $K$ (with respect to a given subfield of $(N-1)$-dimensional constants $E$). 
\enddefinition 

We mention the following simple properties:

a) for any $(K,E)\in\LC $, 
there is always a closed subfield $K_0$ in $K$ containing $E$ 
such that $(K_0,E)\in\LC $ is standard; this field $K_0$ appears in the form 
$E\{\{t\}\}$ with a suitably chosen element $t$ of $\Cal O_K$;

b) if $(\widetilde{K},E)\in\LC (N)$ and $K$ is a closed subfield in $\widetilde{K}$ such that 
$K\supset E$ and $(K,E)\in\LC (N)$, then $(K,E)$ is standard;  
(One can see easily, that $[\widetilde{K}:K]<\infty $ and if $\tilde t_N$ is an $N$-th local parameter for 
$\widetilde{K}$ then $N_{\widetilde{K}/K}\tilde t_N$ is an $N$-th local 
parameter for $K$.)

c) if $(K,E)\in\LC $ is standard then for any finite extension $E'$ of $E$, 
$(KE',E')\in\LC $ is standard; (Any $N$-th local parameter in $K$ is still an 
$N$-th local parameter in $KE'$.)

d) any $(K,E)\in\LC (1)$ is standard;

e) for any $(K,E)\in\LC (2)$, there is a finite extension $E'$ of $E$ 
such that $(KE',E')\in\LC (2)$ is standard. (This follows from 
Epp's Theorem.)
\medskip 

The following property plays a very important role in the construction of  
ramification theory for higher dimensional fields. 

\proclaim{Proposition 2.1} Suppose $(K,E), (L,E)\in\LC (N)$, 
$L\supset K$ and $(L,E)$ is standard. Then 
$\Cal O_L=\Cal O_K[t_N]$, where $t_N$ is an $N$-th 
local parameter in $L$. 
\endproclaim 

\demo{Proof} Clearly, $\Cal O_K[t_N]\subset\Cal O_L$. 

Let $t_1,\dots ,t_{N-1}$ be local parameters in $E$. It will be sufficient 
to prove that 
$$t_1^{a_1}\dots t_{N-1}^{a_{N-1}}t_N^{a_N}\in\Cal O_K[t_N]$$
if $(a_1,\dots ,a_{N-1},a_N)\geqslant \bar 0_N$. 

We can assume that $a_N<0$ (otherwise, there is nothing to prove).

Notice that $\tilde t_N=N_{L/K}t_N$ is an $N$-th 
local parameter for $K$ and 
$\tilde t_Nt_N^{-1}\in\Cal O_K[t_N]$. Therefore, 

$$t_1^{a_1}\dots t_{N-1}^{a_{N-1}}t_N^{a_N}=
t_1^{a_1}\dots t_{N-1}^{a_{N-1}}\tilde t_N^{a_N}
(\tilde t_Nt_N^{-1})^{-a_N}\in\Cal O_K[t_N]$$
because $t_1^{a_1}\dots t_{N-1}^{a_{N-1}}\tilde t_N^{a_N}\in\Cal O_K$. 
The proposition is proved. 
\enddemo

2.2. The following theorem plays in our setting a role of a higher dimensional 
version of Epp's Theorem.

\proclaim{Theorem 1} If $(K,E)\in\LC (N)$, then there is a finite 
separable extension 
$E'$ of $E$ such that $(KE',E')\in\LC (N)$ is standard.
\endproclaim 

\demo{Proof} Use induction on $N$.   

If $N=1$ there is nothing to prove. Notice that the case $N=2$ follows from 
Epp's Theorem. 

 Suppose $N>1$ and the theorem holds for local fields of 
dimension $<N$.

\proclaim{Proposition 2.2} Suppose $(K,E)\in\LC (N)$, 
then 
there is a finite separable 
extension $\widetilde{E}$ of $E$ such that if $\widetilde{K}=K\widetilde{E}$ then 
\newline 
{\rm 1)} $\widetilde{K}$ and $\widetilde{E}$ have a common first uniformiser;
\newline 
{\rm 2)} $\widetilde{E}^{(1)}$ is algebraically closed in $\widetilde{K}^{(1)}$.
\endproclaim 

This proposition will be proved in nn.2.3-2.9 below. It implies the statement of 
Theorem 1 as follows.

By the above property 2), $(\widetilde{K}^{(1)},\widetilde{E}^{(1)})\in\LC (N-1)$. Therefore, 
by the inductive assumption 
there is a finite separable extension $E_1$ of $\widetilde{E}^{(1)}$ 
such that $(K_1,E_1)$ is standard (where $K_1=\widetilde{K}^{(1)}E_1$). 
Denote by $\bar t_2,\dots ,\bar t_{N}$ a system of local parameters in 
$K_1$ such that $\bar t_2,\dots ,\bar t_{N-1}$ is a system of local parameters of 
$E_1$. 
Let $E'$ be an unramified extension $\widetilde{E}$ such that ${E'}^{(1)}=E_1$. 
Notice that if $K'=KE'$ then ${K'}^{(1)}=K_1$. Let $t_2,\dots ,t_{N-1}$ 
be liftings of $\bar t_2,\dots ,\bar t_{N-1}$ to $O^{(1)}_{E'}$ and let 
$t_N$ be a lifting of $\bar t_N$ to $O^{(1)}_{K'}$. Then 
$t_1,\dots ,t_N$ is a system of local parameters in $K'$ and 
$t_1,\dots ,t_{N-1}$ is a system of 
local parameters in $E'$, i.e. $(K',E')\in\LC (N)$ is standard.  
\enddemo

\demo{Proof of Proposition 2.2}
\enddemo 

2.3. Choose a standard $(K_0,E)\in\LC (N)$ such that $K_0\subset K$, and denote by 
$t_1,\dots ,t_N$ a system of local parameters in $K_0$ such that the 
first $N-1$ of them give a system of local parameters in $E$. 

It will be sufficient to prove our theorem 
for 
extensions $K/K_0$ satisfying one of the following 
conditions  (because any finite extension 
of $K_0$ can be embedded into a bigger extension obtained as a sequence 
of such subextensions): 
\medskip 

$a_0$) there is a finite extension $\widetilde{E}$ of $E$ such that  
$\widetilde{K}:=K\widetilde{E}$ is unramified over 
$\widetilde{K}_0:=K_0\widetilde{E}$, i.e. such that both fields $\widetilde{K}$ and $\widetilde{K}_0$ 
have the same 
first uniformiser and $\widetilde{K}^{(1)}$ is separable over $\widetilde{K}_0^{(1)}$; 
\medskip 

$a_1$) $K/K_0$ is a cyclic extension of a prime to $p$ degree $m$;
\medskip

b) $K/K_0$ is a cyclic extension of degree $p$ such that after arbitrary 
finite extension of $E$ the corresponding extension of 
first residue fields is either trivial or purely inseparable. 
When considering this case below we shall treat separately 
the subcases: 

\ \ \ $b_1$) $\char K=0$;

\ \ \ $b_2$) $\char K=p$. 
\medskip 

c) $K/K_0$ is a purely non-separable extension of degree $p$. 
\medskip 

Following the terminology from [Zh2] we can call $(K,E)$ an almost constant 
extension of $(K_0,E)$ in the case a) and an infernal elementary extension 
in the case b). 
\medskip 

\subsubhead{\rm 2.4} The case $a_0)$ 
\endsubsubhead

This case easily follows from the following observation.  
Consider the natural field embedding 
$\widetilde{E}^{(1)}\subset \widetilde{K}_0^{(1)}$. Then 
$(\widetilde{K}_0^{(1)},\widetilde{E}^{(1)})\in\LC (N-1)$. Indeed, 
$(\widetilde{K}_0,\widetilde{E})\in\LC (N)$ is standard, 
then $\widetilde{E}^{(1)}$ is a field of $(N-2)$-dimensional 
constants in $\widetilde{K}_0^{(1)}$, which is  
algebraically closed 
in $\widetilde{K}_0^{(1)}$. On the other hand, 
$\widetilde{E}^{(1)}$ is separably closed in $\widetilde{K}^{(1)}$ 
(otherwise, $\widetilde{E}$ will have a non-trivial unramified extension in $\widetilde{K}$). 
This implies that any finite extension $E'$ of $\widetilde{E}^{(1)}$ in 
$\widetilde{K}^{(1)}$ is either purely inseparable or trivial. 
Therefore, 
$E'\subset\widetilde{K}_0^{(1)}$ (because $\widetilde{K}^{(1)}/\widetilde{K}_0^{(1)}$ 
is separable) and $E'=\widetilde{E}^{(1)}$ 
(because $\widetilde{E}^{(1)}$ is algebraically closed in $\widetilde{K}^{(1)}_0$). 
\medskip 

\subsubhead{\rm 2.5.} The case $a_1)$ 
\endsubsubhead   

We can assume that $E$ contains a primitive $m$-th root of unity. Then 
$K=K_0(\root m \of {t_1^{a_1}\dots t_N^{a_N}})$, where 
$a_1,\dots ,a_N\in\Bbb Z_{\geqslant 0}$, and we can assume that $\operatorname{gcd}(a_N,m)=1$. 
Let $\widetilde{E}=E(\root m\of {t_1},\dots ,\root m\of {t_{N-1}})$, then 
$\widetilde{E}$ has local parameters $\root m\of {t_1},\dots ,\root m\of {t_{N-1}}$ and this system can 
be extended to a system of local parameters in $\widetilde{K}=K\widetilde{E}$ 
by adding $\root m\of {t_N}$. So, $(\widetilde{K},\widetilde{E})$ is 
standard and $\widetilde{E}^{(1)}$ 
is algebraically closed in $\widetilde{K}^{(1)}$. 
\medskip 

\subsubhead{\rm 2.6} Special extensions 
\endsubsubhead 

For our future targets we need to keep control on the choice of 
the extension $\widetilde{E}$ of $E$ in the proposition 2.2. This idea goes back to 
the paper [ZhK] where it was proved that Epp's elimination of wild ramification for 
an infernal extension can be done by the use of subextensions of a 
given deeply ramified extension.

Consider an increasing sequence of finite extensions 
$$E\subset \widetilde{E}_0\subset E_0\subset 
\widetilde{E}_1\subset E_1\subset \dots 
\subset\widetilde{E}_n\subset E_n\subset\dots $$
such that each $\widetilde{E}_n$ and $E_n$ have a system of local parameters 
$\tilde t_{1n},\dots ,\tilde t_{N-1,n}$ and, respectively, 
$t_{1n},\dots ,t_{N-1,n}$, satisfying the following condition:  
\medskip 

{\bf C.}\ \ {\it 
There is a $c>0$ such that for all 
$1\leqslant i\leqslant N-1$ and $n\geqslant 1$, 
$$v^1\left (\frac{t_{in}^p}{\tilde t_{i,n-1}}-1\right )\geqslant c$$ 
where $v^1$ is a $t_1$-adic (1-dimensional) 
valuation on $\bar K$ normalised by the 
condition $v^1(t_1)=1$.}
\medskip 

Proposition 2.2 will be implied in the cases $b)$ and $c)$ by the following statement. 

\proclaim{Proposition 2.3} Suppose that $K, K_0$ and $E$ satisfy the 
assumptions from the cases $b)$ or $c)$. Then there is 
an $n^*\in\Bbb Z_{\geqslant 0}$ (depending only on 
the extension $K/K_0$ and the $c$ from the above condition C) 
such that proposition 2.1 
holds with $E'=E_{n^*}$. 
\endproclaim 

\subsubhead {\rm 2.7.} The case $b_2)$  
\endsubsubhead 
 
In the case $b_2)$ we have    
$K=K_0(\theta )$, $\theta ^p-\theta =\xi $, where   
$\xi\in K_0$ is a power series 
$$\xi =\sum\Sb \bar a \endSb [\alpha _{\bar a}]t_1^{a_1}\dots t_N^{a_N}$$
with restrictions on its coefficients 
described in the beginning of section 1. Applying the Artin-Schreier equivalence we 
can assume also that 
it contains only non-zero terms with $\bar a\leqslant \bar 0_N$ 
and $\bar a\not\equiv 0\operatorname{mod}p$ if $\bar a\ne 0$.

Set $\xi =\xi '+\xi ''$, where 
$$\xi '=\sum\Sb a_N=0 \endSb [\alpha _{\bar a}]t_1^{a_1}\dots t_{N-1}^{a_{N-1}}\ \ \ \ 
\xi ''=\sum\Sb a_N\ne 0 \endSb [\alpha _{\bar a}]t_1^{a_1}\dots t_{N-1}^{a_{N-1}}t_N^{a_N}$$

Let 
$$A=\operatorname{min}\{a_1\ |\ \alpha _{\bar a}\ne 0, a_N=0\}=v^1(\xi ')$$
$$B=\operatorname{min}\{a_1\ |\ \alpha _{\bar a}\ne 0, a_N\ne 0\}=v^1(\xi ''),$$ 
where 
 $v^1$ is a $t_1$-adic valuation from the above condition $C$.

Notice that the first set can be empty. In this case we set by definition 
$A=0$. The second set is never empty: otherwise, $K$ is a composite 
of an algebraic extension of $E$ and $K_0$, i.e. $E$ is not 
algebraically closed in $K$. For any $s\in\Bbb Z_{\geqslant 0}$, 
let 
$$B^{(s)}=\operatorname{min}\{a_1\ |\ \alpha _{\bar a}\ne 0, v_p(a_N)=s\}$$ 
(we set $B^{(s)}=0$ if the corresponding subset of indices is empty). 
Then $B=\operatorname{min}\{B^{(s)}\ |\ s\geqslant 0\}$.

Notice that if we pass from 
$E$ to its finite extension $\widetilde{E}_0$, cf. condition C, then 
$\tilde t_{10},\dots ,\tilde t_{N-1,0},t_N$ is a system of local 
parameters for $K_0\widetilde{E}_0$. Rewrite $\xi $ in terms of these 
local parameters and apply to this expression the Artin-Schreier equivalence  
to get rid of all $p$-th powers and terms from the maximal ideal of 
$O_{K_0\widetilde{E}_0}$. This procedure gives an analogue 
$\tilde\xi _0$ of $\xi $ for the extension $K\widetilde{E}_0/K_0\widetilde{E}_0$. 
As earlier, use the $t_1$-adic valuation $v^1$ to define the analogues 
$\widetilde{A}_0$, $\widetilde{B}_0$, $\widetilde{B}_0^{(s)}$ of, respectively, 
$A$, $B$ and $B^{(s)}$, $s\geqslant 0$. 

\proclaim{Lemma 2.4} {\rm a)} $\widetilde{A}_0\geqslant A$;
\newline 
{\rm b)} for all $s\geqslant 0$, 
$\widetilde{B}_0^{(s)}\geqslant\operatorname{min}\left\{\frac{1}{p^u}B^{(s+u)}\ |\ u\geqslant 0\right\}$.
\endproclaim 

Apply the similar procedure to the extensions $E_0,\widetilde{E}_1,E_1,\dots $ to get the 
invariants 
\linebreak $A_0,B_0,B_0^{(s)}$, 
$\widetilde{A}_1$, $\widetilde{B}_1$, $\widetilde{B}_1^{(s)}$,  $A_1$, $B_1$, $B_1^{(s)}$,\dots . 

Similarly, we have the following property.

\proclaim{Lemma 2.5} For all $i,s\geqslant 0$, 
\newline 
{\rm a)} $\widetilde{A}_{i+1}\geqslant A_i$;
\newline 
{\rm b)} $\widetilde{B}_{i+1}^{(s)}\geqslant\operatorname{min}
\left\{\frac{1}{p^u}B_i^{(s+u)}\ |\ u\geqslant 0\right\}$.
\endproclaim 

When passing through the special extensions $E_i/\widetilde{E}_i$, $i\geqslant 0$, 
we have the better estimates:

\proclaim{Lemma 2.6} For all $i\geqslant 0$ and $s\geqslant 1$, 
\newline 
{\rm a)} $A_i\geqslant\operatorname{min}\left\{\frac{1}{p}\widetilde{A}_i, \widetilde{A}_i+c\right\}$;
\newline 
{\rm b)} $B_i^{(0)}\geqslant \operatorname{min}\left\{ \widetilde{B}_i^{(0)}; \frac{1}{p}\widetilde{B}_i^{(1)}; 
\frac{1}{p^u}(\widetilde{B}_i^{(u)}+c), u\geqslant 2\right\}$;
\newline 
{\rm c)} $B_i^{(s)}\geqslant\operatorname{min}\left\{\frac{1}{p}\widetilde{B}_i^{(s+1)}; 
\frac{1}{p^u}\left (\widetilde{B}_i^{(s+u)}+c\right ), u\geqslant 0\right\}$.
\endproclaim 

\proclaim{Corollary 2.7} {\rm a)} $\lim_{i\to\infty }A_i=0$;
\newline 
{\rm b)} if $\gamma _i=\operatorname{min}\{ B_i^{(s)}\ |\ s\geqslant 1\}$, 
then $\lim_{i\to\infty }\gamma _i=0$.
\endproclaim 

\proclaim{Lemma 2.8} If $i\geqslant 0$ is such that $B_i^{(0)}<B_i^{(s)}$ for all $s\geqslant 1$, 
then for all $u\geqslant i$, $B_u=B_i^{(0)}$.
\endproclaim 

\proclaim{Corollary 2.9} There is an index $n^*$ such that $A_{n^*}>B_{n^*}$. 
\endproclaim

 So, if $n\geqslant n^*$,  then $K_{n}^{(1)}=K_{0n}^{(1)}(\bar\theta )$ 
with 
$$\bar\theta ^p=\bar\eta :=(t_{1n}^{-Be_n}\xi _n)\operatorname{mod}m_{K_{0n}}^{(1)}$$
where $B=B_n=B_{n^*}$ and $e_n^{-1}=v^1(t_{1n})$. 
Clearly, $\bar\eta\notin {K_{0n}^{(1)}}^p+E_{n}^{(1)}$ and, therefore, 
$E_{n}^{(1)}$ remains to be algebraically closed in $K_{n}^{(1)}$. 

Besides, if $n\geqslant \operatorname{min}\{n^*,1\}$ then the first uniformiser 
$t_{1n}$ appears in the leading term of the 
$\xi _{n}''$ with an exponent divisible by $p$ and, therefore, it is also a   
uniformiser for $K_{n}$. So, proposition 1.3 is  proved in the case $b_2)$. 
\medskip  

\subsubhead {\rm 2.8} The case $c)$ 
\endsubsubhead 

In this case we have    
$K=K_0(\theta )$, $\theta ^p=\xi $, where   
$\xi\in K_0$ is the power series from n.2.7, 
containing non-zero terms only with $\bar a\not\equiv 0\operatorname{mod}p$.

Set $\xi =\xi '+\xi ''$, where 
$$\xi '=\sum\Sb a_N\equiv 0\operatorname{mod}p \endSb 
[\alpha _{\bar a}]t_1^{a_1}\dots t_{N-1}^{a_{N-1}}t_N^{a_N}\ \ \ \ 
\xi ''=\sum\Sb a_N\not\equiv  0\operatorname{mod}p \endSb 
[\alpha _{\bar a}]t_1^{a_1}\dots t_{N-1}^{a_{N-1}}t_N^{a_N}$$

Let 
$$A=\operatorname{min}\{a_1\ |\ \alpha _{\bar a}\ne 0, a_N\equiv 0\operatorname{mod}p\}=v^1(\xi ')$$
$$B=\operatorname{min}\{a_1\ |\ \alpha _{\bar a}\ne 0, a_N\not\equiv 0\operatorname{mod}p\}=v^1(\xi ''),$$ 
where 
 $v^1$ is a $t_1$-adic valuation from the above condition $C$.

Notice that the first set can be empty. In this case we set by definition 
$A=+\infty $. The second set is never empty: otherwise, $\theta $ is algebraic over $E$, 
i.e. $E$ is not 
algebraically closed in $K$.

If we pass from 
$E$ to its finite extension $\widetilde{E}_i$, where 
$i=0,1,\dots $, cf. condition C, then 
$\tilde t_{i1},\dots ,\tilde t_{i,N-1},t_N$  
is a system of local 
parameters for $K_0\widetilde{E}_i$. Rewrite $\xi $ in terms of these 
local parameters and take away   
all $p$-th power terms. This procedure gives an analogue 
$\tilde\xi _i$ of $\xi $ for the extension $K\widetilde{E}_i/K_0\widetilde{E}_i$. 
As earlier, use $t_1$-adic valuation $v^1$ to define the analogues 
$\widetilde{A}_i$ and $\widetilde{B}_i$.

Similarly, introduce the invariants $A_i$ and $B_i$, where $i=0,1,\dots $, 
when passing in the above procedure from $E$ to $E_i$.  

We have the following estimates. 

\proclaim{Lemma 2.10} {\rm a)} $\widetilde{A}_0\geqslant A$ and $\widetilde{B}_0=B$;
\newline 
{\rm b)} for all $i\geqslant 0$, $\widetilde{A}_{i+1}\geqslant A_i$ and 
$\widetilde{B}_{i+1}=B_i$;
\newline 
{\rm c)} for all $i\geqslant 0$,  $A_i\geqslant \widetilde{A}_i+c$ and $B_i=\widetilde{B}_i$.
\endproclaim 

This implies immediately that there is an index $n^*$ such that $A_{n^*}>B_{n^*}=B$. 
Therefore, for all $n\geqslant n^*$, $K_n^{(1)}=K_{0n}^{(1)}(\bar\theta )$, 
where $\bar\theta ^p$ is the image of $t_{1n}^{-Be_n}\xi _n$ in $K_{0n}^{(1)}$, where 
$e^{-1}_n=v^1(t_{1n})$. 
Clearly, $E_n^{(1)}$ is still algebraically closed in $K_n^{(1)}$.  
Even more, 
if 
$n\geqslant\min\{n^*,1\}$ then $t_{1n},\dots ,t_{N-1,n}$ appear in the leading term of 
$\xi _n''$ with divisible by $p$ exponents. In particular, $t_{1n}$ is still a (first) uniformiser for 
$K_n$. 

The case $c)$ is also considered. 
\medskip

\subsubhead{\rm 2.9.} Characteristic 0 analogue of the Artin-Schreier theory 
\endsubsubhead 

 The characteristic 0 case $b_1)$ can be treated similarly to the characteristic 
$p$ case $b_2)$ due to the characteristic 0 analogue of the Artin-Schreier theory 
from [Ab1]. This construction can be briefly reminded as follows. 

Suppose $L_0$ is a complete discrete valuation field of 
characteristic 0 with the maximal ideal $\m _{L_0}$ and the 
residue field $k$ of characteristic $p$. Assume that 
$\zeta _p\in L_0$ (where $\zeta _p$ is a primitive $p$-th root of unity) 
and let $\pi _1\in L_0$ be such that $\pi _1^{p-1}=-p$. 

\proclaim{Proposition 2.11} \ \ 
\newline 
{\rm a)} $L=L_0(\root p\of v)$ with $v\in 1+\pi _1\m _{L_0}$ if and only if 
$L=L_0(\theta )$, where $\theta ^p-\theta =w$ with $w\in p^{-1}\m _{L_0}$;
\newline 
{\rm b)} With the above notation and assumptions $L$ admits another presentation 
$L=L_0(\theta _1)$, 
where $\theta _1^p-\theta _1=w_1\in p^{-1}\m _{L_0}$, if 
$w_1=w+\eta ^p-\eta $ with $\eta\in L_0$ such that $\eta ^p\in p^{-1}\m _{L_0}$.
\endproclaim 

\demo{Proof} We only sketch the idea of the proof. 

Let 
$$E(X)=\exp \left (X+X^p/p+\dots +X^{p^n}/p^n+\dots \right )\in\Bbb Z_p[[X]]$$ 
be the Artin-Hasse exponential. 
Then $v=E(\pi _1V)$ with $V\in \m _{L_0}$ and if $u^p=v$, $u\in L$, then $u=E(U)$ with 
$U\in \m _L$. Then the equivalence 
$$E(X^p)=E(X^p)\exp (pX)\equiv E(X^p+pX)\operatorname{mod}(p^2X,pX^p)$$
implies that 
$$U^p+pU \equiv \pi _1V\operatorname{mod}\pi _1p\m _L$$
(notice that $U^p\in \pi _1\m _L$). 

Divide both 
sides of the above equivalence by $\pi _1^{p}$ and deduce that 
$L=L_0(\theta )$, where $\theta ^p-\theta =w\in p^{-1}\m _{L_0} $ with 
$\theta\equiv \pi _1^{-1}U\operatorname{mod}\m _{L}$ and 
$w\equiv p^{-1}V\operatorname{mod}\m _{L_0}$. 
\enddemo 
\medskip 

\subsubhead{\rm 2.10} The case $b_1)$
\endsubsubhead

\proclaim{Proposition 2.12} 
Suppose $K, K_0$ and $E$ satisfy the condition $b_1)$ from 
n.2.3. Then there is an $n^*\in\Bbb Z_{\geqslant 0}$ such that 
proposition 2.3 holds with $E'=E_{n^*}$. 
\endproclaim

\demo{Proof} Assume first that $\zeta _p\in E$. 

Then $K\widetilde{E}_0=(K_0\widetilde{E}_0)(\root p\of {\tilde v_0})$, where 
$\tilde v_0=\tilde t_{10}^{c_1}\dots \tilde t_{N0}^{c_N}(1+\tilde a)$, 
$\tilde a\in \m _{K\widetilde{E}_0}$ and $c_1,\dots ,c_N\in\Bbb Z_{\geqslant 0}$. 

Then the condition $C$ from n.2.6 implies that 
$$\tilde v_0=t_{10}^{pc_1}\dots t_{N0}^{pc_N}(1+\tilde a^p+\tilde b),$$
where $\tilde a,\tilde b\in \m _{KE_0}$ and $v^1(\tilde b)\geqslant c$. This implies that 
$KE_0=K_0E_0(\root p\of {v_0})$, where 
$v_0=1+pa+b$ with $a,b\in \m _{K_0E_0}$ such that $v^1(b)\geqslant c$. 

By continuiing the above procedure we  obtain that $KE_n=(K_0E_n)(\root p\of {v_n})$ where 
$v_n=1+pb_n$ with $b_n\in \m _{K_0E_n}$. Since $v_n\in 1+\pi _1\m _{K_0E_n}$, the extension 
$KE_n/K_0E_n$ can be given via the analogue of the Artin-Schreier theory from 
n.2.9 and we can proceed further as in n.2.7 to finish the proof of our proposition. 

Suppose now that $\zeta _p\notin E$. 

Let $K'=K(\zeta _p)$, $K_0'=K_0(\zeta _p)$, 
$E'=E(\zeta _p)$ and $\widetilde{E}_n'=\widetilde{E}_n(\zeta _p)$, 
$E_n'=E_n(\zeta _p)$ for all $n\geqslant 0$. Then the tower 
$$\widetilde{E}'\subset \widetilde{E}_0'\subset E_0'\subset\dots \widetilde{E}_n'\subset E_n'\subset\dots $$
satisfies the condition $C$ from n.2.6. Therefore, there is an $n^*$ such that 
if $K_n'=KE_n'$ and $K'_{n0}=K_0E'_n$, then 
${E_n'}^{(1)}$ is algebraically closed in ${K_n'}^{(1)}$. 

Let $F$ be a non-trivial purely inseparable extension of $E_n^{(1)}$ in $K_n^{(1)}$. 
Then $F{E'_n}^{(1)}$ is a non-trivial purely inseparable extension of 
${E'_n}^{(1)}$ in ${K_n'}^{(1)}$ 
\linebreak 
(use that $[{E'_n}^{(1)}:E_n^{(1)}]<p$). But this contradicts 
to the fact that  ${E_n'}^{(1)}$ is 
algebraically closed in ${K_n'}^{(1)}$. Therefore, 
$E_n^{(1)}$ is algebraically closed in $K_n^{(1)}$, because 
it is its separably closed subfield. 

The proposition is completely proved.
\enddemo 
\medskip

\subhead 3. Ramification theory and Krasner's Lemma 
\endsubhead 
\medskip 

\subsubhead {\rm 3.1} The category of local fields with $F$-structures 
\endsubsubhead 
\medskip 

This category $\LF (N)$ will apear as the disjoint union of its two full 
subcategories $\LF _0(N)$ and $\LF _p(N)$. 
\medskip 

{\it The category $\LF _0(N)$.} 

Choose a simplest $N$-dimensional local field 
$L_0=\Bbb Q_p\{\{t_N\}\}\dots \{\{t_2\}\}$. Define its $F$-structure as an increasing  
sequence of closed subfields  
$\{L_{0}(i)\ |\ 1\leqslant i\leqslant N\}$ 
with the system of local parameters $p=t_1,t_2, \dots ,t_N$. 
Choose an algebraic closure $\bar L_0$ of $L_0$. Denote by $\Bbb C(N)_p$ 
the completion of $\bar L_0$ with respect to its first ($p$-adic) valuation. 
For $1\leqslant i\leqslant N$, denote by $\Bbb C(i)_p$ the completion 
of the algebraic closure of $L_{0}(i)$ in $\Bbb C(N)_p$. 
It will be convenient to have a special agreement for $i=0$. 
By definition, $\Bbb C(0)_p$ is the completion of the maximal unramified 
extension of $\Bbb Q_p$ in $\Bbb C(N)_p$ and $L_0(0)=L_0\cap \Bbb C(0)=\Bbb Q_p$.  
Notice that $\Bbb C(1)_p=\Bbb C_p$ is the $p$-adic completion of 
an algebraic closure of $\Bbb Q_p$.

Clearly, the $P$-topological structure of finite  extensions of $L_0$ 
induces the $P$-topological structures on the fields 
$\Bbb C(0)_p\subset\Bbb C(1)_p\subset\dots \subset\Bbb C(N)_p$.

The objects of the category $\LF _0(N)$ are 
finite extensions $K$ of $L_0$ in $\Bbb C(N)_p$ with the induced 
$F$-structure. This structure is given by the sequence of 
algebraically closed and $P$-closed subfields 
$\{K(i)\ |\ 0\leqslant i\leqslant N\}$, where 
$K(i)=K\cap\Bbb C(i)_p$. Notice that $K(0)$ is the maximal 
unramified extension of $\Bbb Q_p$ in $K$. We agree 
to use the notation $\bar K$ for the algebraic closure of $K$ 
in $\Bbb C(N)_p$. Notice that $\Gamma _K=\Aut (\bar K/K)$ consists of $P$-continuous 
field automorphisms $\tau $ of $\Bbb C(N)_p$ such that $\tau |_K=\id $ and for 
all $0\leqslant i\leqslant N$, $\tau (\Bbb C(i)_p)=\Bbb C(i)_p$. It is well-known [Hy], 
that  $\Bbb C(N)_p^{\Gamma _K}=K$ and, therefore, 
for all $0\leqslant i\leqslant N$, 
$\Bbb C(i)_p^{\Gamma _K}=K(i)$.

Suppose $K,L\in\LF _0(N)$. Then the corresponding 
set of morphisms 
\linebreak 
$\operatorname{Hom}_{\LF (N)}(K,L)$ 
consists of all 
$P$-continuous field morphisms $\varphi :\Bbb C(N)_p\rightarrow\Bbb C(N)_p$ 
such that  for $0\leqslant i\leqslant N$, 
\newline 
{\rm a)} $\varphi (\Bbb C(i)_p)=\Bbb C(i)_p$;
\newline 
{\rm b)} $\varphi (K)\subset L$.  

Notice that any $\varphi\in\operatorname{Hom}_{\LF (N)}(K,L)$ 
transforms the $F$-structure of $K$ to the $F$-structure of $L$. 
\medskip 

{\it The category $\LF _p(N)$}. 

We proceed similarly to the above characteristic 0 case. 
Choose a basic $N$-dimensional local field 
$L_p=\Bbb F_p((t_N))\dots ((t_1))$ and define its  
$F$-structure by a sequence of subfields $\{L_p(i)\ |\ 0\leqslant i\leqslant N\}$ 
such that $L_p(i)$ has local parameters $t_1,\dots ,t_i$. 
Choose an algebraic closure $\bar L_p$ of $L_p$. Denote by $\Cal C(N)_p$ 
the completion of $\bar L_p$ with respect to its first valuation. 
For $0\leqslant i\leqslant N$, denote by $\Cal C(i)_p$ the completion 
of the algebraic closure of $L_p(i)$ in $\Cal C(N)_p$. 
As earlier, the $P$-topological structure of finite  extensions of $L_p$ 
induces the $P$-topological structures on the fields 
$\bar\Bbb F_p=\Cal C(0)_p\subset\Cal C(1)_p\subset\dots \subset\Cal C(N)_p$.

The objects of the category $\LF _p(N)$ are 
finite extensions $K$ of $L_p$ in $\Cal C(N)_p$ 
with the induced $F$-structure 
$\{K(i)\ |\ 0\leqslant i\leqslant N\}$, where 
$K(i)=K\cap\Cal C(i)_p$. Notice that $\Cal C(N)_p^{\Gamma _K}=\Cal R(K)$ 
--- the radical closure (=the completion of the maximal 
purely non-separable extension) of $K$ in $\Cal C(N)_p$. 
Similarly for $0\leqslant i\leqslant N$, it holds that 
$\Cal C (i)_p^{\Gamma _K}=\Cal R(K(i))$. The morphisms in 
$\LF _p(N)$ are defined also along lines in the above charactersitic 0 case. 
\medskip

\subsubhead{\rm 3.2} Standard $F$-structure 
\endsubsubhead 

We say that the $F$-structure on $L\in\LF (N)$ is standard if there is a system 
of local parameters $t_1,\dots ,t_N$ in $L$ such that 
for all $1\leqslant r\leqslant N$, $t_1,\dots ,t_r$ is  
a system of local parameters for $L(r)$. Applying the above Theorem 1  
we obtain easily  the following 

\proclaim{Proposition 3.1} 
For any $E\in\LF (N)$,  there is a finite separable extension $E'$ of 
$E(N-1)$ such that $EE'$ has a standard $F$-structure.  
\endproclaim 

\remark{Remark} 
The above proposition played a fundamental role in the construction 
of the higher dimensional ramification theory in [Ab5], but its proof in [Ab5] was not complete, due 
to reasons mentioned in the Remark from 2.1.  Notice that the construction of ramification 
theory, cf. n.3.3 below, can be based only on the result of Theorem 1. 
\endremark

Note that the 
$F$-structure allows to treat higher dimensional local fields 
in a very similar way to classical complete discrete valuation fields 
with finite residue fields. For example, for any finite extension 
of local fields with $F$-structure we can introduce: 
\medskip  

a) {\it a vector ramification index} $\bar e(L/K)=(e_1,\dots ,e_N)$.  

Any finite extension of $K$ in $\bar K$ 
appears with a natural $F$-structure and a natural $P$-toplogy. 
In particular, if $L\subset M$ are such subfields in $\bar K$ then 
its vector ramification index equals 
$e(M/K)=(e_1,\dots ,e_N)$, where for $1\leqslant r\leqslant N$, 
$e_r=[M(r):L(r)]/[M(r-1):L(r-1)]$. This index plays a role of the usual 
ramification index in the theory of 1-dimensional local fields. 
\medskip

b) {\it a canonical $N$-valuation }$v_{L}:L\longrightarrow\Bbb Q^N\cup\{\infty \}$.

If $L$ has a standard $F$-structure and 
$t_1,\dots ,t_N$ is a corresponding system of local parameters, then $v_L$ 
is uniquely defined by the conditions 
$v_L(t_1)=(1,0,\dots ,0)$, $v_L(t_2)=(0,1,0,\dots ,0)$,\dots , 
$v_L(t_N)=(0,0,\dots ,0,1)$. Otherwise, one should use a finite extension  
$L_1$ of $L$ with standard $F$-structure and set $v_L=\bar e(L_1/L)^{-1}v_{L_1}$. 
\medskip 

\subsubhead {\rm 3.3}
Review of ramification theory  
\endsubsubhead 

Suppose $K\in\LF (N)$. 
   Then $\Gamma _K=\Aut (\bar K/K)$ has a canonical decreasing 
filtration by ramification subgroups  
$\{\Gamma _K^{(j)}\ |\ j\in J(N)\}$ with  
the set of indices
$J(N)=\mathbin{\underset{1\leqslant r\leqslant N}\to\coprod}J_r$. Here 
$J_r=\{ j\in \Bbb Q^r \ |\  j\geqslant \bar 0_r\}$  
with respect to the lexicographic ordering on $\Bbb Q^r$, where 
$\bar 0_r=(0,\dots ,0)\in\Bbb Q^r$. By definition, if $r_1>r_2$ then 
any element from 
$J_{r_1}$ is bigger than any element from 
$J_{r_2}$.

The definition of this filtration can be described as follows.

Let $E/K$ be a finite extension in $\bar K$ 
(this is a subfield in $\Bbb C(N)_p$ or $\Cal C(N)_p$). Consider the finite set   
$I_{E/K}$ of all $P$-continuous embeddings of 
$E$ into $\bar K$ which are the identity on $K$. 

There is a natural filtration of this set 
$$I_{E/K}\supset I_{E/K,0}\supset I_{E/K,(0,0)}\supset\dots \supset I_{E/K,\bar 0_N}$$
where for $1\leqslant r\leqslant N$, 
$I_{E/K,\bar 0_r}$ are embeddings which are the identity on the 
subfield of $(r-1)$-dimensional constants $E(r-1)$.

For $1\leqslant r\leqslant N$ and $j\in J_r$,  
define the set $I_{E/K,j}\subset I_{E/K,\bar 0_r}$ as follows.

Take a suitable finite extension $E'$ of $E(r-1)$ in $\bar K$ such that 
if $\widetilde{E}(r)=E'E(r)$ and $\widetilde{K}(r)=K(r)E'$ then 
$\Cal O_{\widetilde{E}(r)}=\Cal O_{\widetilde{K}(r)}[\theta ]$. 
(Recall, if $L\in\LF (r)$ then 
\linebreak 
$\Cal O_L=\{l\in L\ |\ v_L(l)\geqslant\bar 0_r\}$.)    
Then use the natural identification 
$I_{E/K,\bar 0_r}=I_{\widetilde{E}(r)/\widetilde{K}(r)}$ to define the 
ramification filtration of $I_{E/K}$ in lower numbering  
$$I_{E/K,j}=\{ \tau\in I_{\widetilde{E}(r)/\widetilde{K}(r)}
\ |\ v_{E(r)}(\tau (\theta )-\theta )
\geqslant v_{E(r)}(\theta )+j\}.$$

Introduce an analogue of the Herbrand function 
$\varphi _{E/K}:J(N)\longrightarrow J(N)$ by setting for 
$1\leqslant r\leqslant N$ and $j\in J_r$,  

$$\varphi _{E/K}(j)=\bar e^{-1}_{E(r)/K(r)}\int _{\bar 0_r}^{j}|I_{E/K,j}|dj\in J_r.$$

This gives the upper numbering such that for any $j\in J(N)$, 
$I_{E/K}^{(j)}=I_{E/K,\varphi _{E/K}(j)}$. 
As in the classical situation, if $E_2\supset E_1\supset K$, then 
the natural projection 
\linebreak 
$I_{E_2/K}\longrightarrow I_{E_1/K}$ induces for 
any $j\in J(N)$, an 
epimorphic map from $I_{E_2/K}^{(j)}$ onto $I_{E_1/K}^{(j)}$ and 
$\varprojlim I^{(j)}_{E/K}=\Gamma _K^{(j)}$ is the ramification subgroup of 
$\Gamma _K$ with the upper number $j$.

As an example, consider the case of an extension 
$E/K$ in $\LF (N)$ such that 
$[E:K]=p^N$ and $e(E/K)=(p,\dots ,p)\in\Bbb Q^N$. Then for 
$1\leqslant r\leqslant N$, there are 
$\alpha _r>\bar 0_r$ such that  
for all $j\in J_r$, 
$$\varphi _{E/K}(j)=
\cases j, \text{ if} j<\alpha _r; \\ \alpha _r+\frac{j-\alpha _r}{p}, 
\text{ if }  j\geqslant\alpha _r
\endcases $$

As in the classical case for any finite extension $E/K$, 
the Herbrand function $\varphi _{E/K}:J(N)\longrightarrow J(N)$ 
is a piece-wise linear function with finitely many edge points. 
Define $i(E/K)\in J(N)$ and $j(E/K)\in J(N)$ as the first and the second coordinates 
of the last edge point of the graph of 
$\varphi _{E/K}$. Notice that 
if $1\leqslant r\leqslant N$ and $j\in J_r$ then 
$j\in J_r$ is an edge point iff $\varphi '_-(j)\ne\varphi '_+(j)$, where 
$\varphi '_-(j)$ and $\varphi '_+(j)$ are slopes of $\varphi _{E/K}$ in the left and right 
neighbourhoods of $j$, respectively. (By definition, 
$\varphi '_-(\bar 0_r)=g_{r0}\bar e^{-1}_{E(r)/K(r)}$, where $g_{r0}=[E(r):K(r)E(r-1)]$.)

If $1\leqslant r\leqslant N$ and $j\in J_r$, then $\varphi '_{-}(j)=g_{-}(j)\bar e^{-1}_{E(r)/K(r)}$ and 
$\varphi '_{+}(j)=g_{+}(j)\bar e^{-1}_{E(r)/K(r)}$, where $g_{-}(j)$ and $g_{+}(j)\in \Bbb N$. 
We shall call $g_{-}(j)/g_{+}(j):=\operatorname{mult}_{E/K}(j)$ 
--- the multiplicity of $\varphi _{E/K}$ in $j\in J_r$. We have:
\medskip 

--- $\operatorname{mult}_{E/K}(j)=1$ if and only if $j$ is not an edge point;
\medskip 

--- $\mathbin{\underset{j\in J(N)}\to\prod}\operatorname{mult}_{E/K}(j)=[E:K]$. 
\medskip 
\medskip

\subsubhead{\rm 3.4.} Krasner's lemma 
\endsubsubhead

Suppose 
 $L,K\in\LF (N)$, $L\supset K$, $L(N-1)=K(N-1)$ and $E$ is a 
finite extension of $L(N-1)$ such that 
$(LE,E)\in\LC (N)$ is standard. 
Then 
$\Cal O_{\widetilde{L}}=\Cal O_{\widetilde{K}}[\theta ]$ where 
$\widetilde{L}=LE$, $\widetilde{K}=KE$ and 
$\theta $ is an $N$-th local parameter in $\widetilde{L}$.

Let $F(T)=T^d+a_1T^{d-1}+\dots +a_d\in\Cal O_{\widetilde{K}}[T]$ be the minimal unitary polynomial 
for $\theta $ over $\widetilde{K}$. 
Denote by $\theta _1=\theta ,\theta _2,\dots ,\theta _d\in\bar K$ all roots of $F(T)$.  
Notice that $v_L(\theta _1)=\dots =v_L(\theta _d)=(0,\dots ,0,1)$. 

In this situation the Krasner Lemma can be given by the following proposition.

\proclaim{Proposition 3.2} 
If $\alpha\in\bar K$ is such that 
$v_K(F(\alpha ))=A+(0,\dots ,0,1)$ with 
$A>\bar 0_N$, then 
\newline 
{\rm 1)} 
there is an index $1\leqslant l_0\leqslant d$ such that 
$v_L(\alpha -\theta _{l_0})=a+(0,\dots ,0,1)$, where $\varphi _{L/K}(a)=A$;
\newline 
{\rm b)} if $A>j(L/K)$ then the above index $l_0$ is unique.
\endproclaim

\demo{Proof} Choose an index $l_0$ such that 
$$v_L(\alpha -\theta _{l_0})=
\operatorname{max}\{v_L(\alpha -\theta _{l})\ |\ 1\leqslant l\leqslant d\}.$$
Let $a\in J_N$ be such that $v_L(\alpha -\theta _{l_0})=a+(0,\dots ,0,1)$. 

\proclaim{Lemma 3.3} 
$v_K(F(\alpha ))=\varphi _{L/K}(a)+(0,\dots ,0,1)$.
\endproclaim 

\demo{Proof of lemma} Let    
$i_1<i_2<\dots <i_s$ be the lower indices which correspond to all jumps 
of the ramification filtration on $I_{L/K}$. 
Then for some integers 
$d=g_0>g_1>\dots >g_{s-1}>g_s=1$  and all   
$1\leqslant i\leqslant n$,  
$v_L(\theta -\theta _i)$ takes $g_0-g_1$ times the value 
$i_1+v_L(\theta )$, \dots , $g_{s-1}-g_s$ times 
the value $i_s+v_L(\theta )$. Notice that 
$i_s=i(L/K)$, $\bar e_{L/K}=(1,\dots ,1,d)$ and if 
$i_t\leqslant a<i_{t+1}$ for some $0\leqslant t\leqslant s$ 
(with the agreements $i_0=\bar 0_N$ and  $i_{s+1}=\infty $) then  
$$\varphi _{L/K}(a)=\bar e^{-1}_{L/K}\left (g_0i_1+
\dots +g_{t-1}(i_t-i_{t-1})+g_t(a-i_t)\right ).$$

Clearly, for all $1\leqslant l\leqslant s$, 
$v_L(\alpha -\theta _l)=\min \{v_L(\alpha -\theta _{l_0}),
v_L(\theta _{l_0}-\theta _l)\}$. 
This implies  
 
$$v_L(F(\alpha ))=\sum\Sb 1\leqslant l\leqslant n\endSb v_L(\alpha -\theta _l)$$

$$=(g_0-g_1)(i_1+v_L(\theta ))+\dots +(g_{t-1}-g_t)(i_t+v_L(\theta _{l_0}))+
g_t(a+v_L(\theta _{l_0}))$$

$$=g_0v_L(\theta _{l_0})+g_0i_1+g_1(i_2-i_1)+\dots +
g_{t-1}(i_t-i_{t-1})+g_t(a-i_t)$$

$$=\bar e_{L/K}\left (v_L(\theta _{l_0})+\varphi _{L/K}(a)\right ).$$
The lemma is proved, because $v_K=\bar e_{L/K}^{-1}v_L$.  
\enddemo 

It remains to prove the part 2) of our proposition.  

Suppose  $\theta _{l_1}$ is a root of $F$ with the same property 
$v_L(\alpha -\theta _{l_1})=a+(0,\dots ,0,1)$. 
Then 
$v_L(\theta _{l_1}-\theta _{l_0})\geqslant a+(0,\dots ,0,1)$. 
But if $A>j(L/K)$ then $a>i(L/K)$ and $\theta _{l_1}=\theta _{l_2}$.

The proposition is proved.
\enddemo 

\proclaim{Corollary 3.4} With the above assumption and notation 
$$v_K(D(F))=(1,\dots ,1,d)j(L/K)-i(L/K)+(0,\dots ,0,d-1)$$
where $D(F)$ is the discriminant of $F$.
\endproclaim 

\demo{Proof} Let $\delta (F)=(\theta -\theta _2)\dots (\theta -\theta _d)$ 
be the different of $F$. Then 
$$v_K(D(F))=v_L(\delta (F))=\sum\Sb 2\leqslant i\leqslant d\endSb 
v_L(\theta -\theta _i)$$
$$=(g_0-g_1)(i_1+v_L(\theta ))+\dots +(g_{s-1}-g_s)(i_s+v_L(\theta ))$$
$$=\bar e(L/K)\varphi _{L/K}(i_s)-i_s+(d-1)v_L(\theta )$$
It remains to notice that $\bar e(L/K)=(1,\dots ,1,d)$, 
$i_s=i(L/K)$ and $\varphi _{L/K}(i_s)=j(L/K)$. 
\enddemo 

\proclaim{Corollary 3.5} $j(L/K)\leqslant 2v_K(D(F))$. 
\endproclaim 

\demo{Proof} This follows from the last corollary because 
$i(L/K)\leqslant v_L(\delta (F))=v_K(D(F))$. 
\enddemo 
\medskip

\subhead 4. Families of increasing towers 
\endsubhead 
\medskip 

In this section we work with local fields of characteristic 0 from $\LF _0(N)$.  
\medskip 

\subsubhead {\rm 4.1.} The category $\Cal B(N)$ 
\endsubsubhead

The objects of $\Cal B(N)$ are  
increasing sequences $K_{\centerdot }=
\{K_n\ |\ n \geqslant 0\}$ of $K_n\in\LF _0(N)$. 
If 
$K_{\centerdot }, L_{\centerdot}\in\Cal B(N)$, then 
$\Hom _{\Cal B(N)}(K_{\centerdot },L_{\centerdot })$ 
consists of field automorphisms 
$f:\Bbb C(N)_p\longrightarrow\Bbb C(N)_p$ 
such that 

--- $f$ is $P$-continuous;

--- $f$ is compatible with $F$-structure;

--- $f(K_n)\subset L_n$ for all $n\gg 0$.

Clearly, if $K_{\centerdot }=\{K_n\ |\ n\geqslant 0\}\in\Cal B(N)$ then for 
any $1\leqslant r\leqslant N$, the subfields of constants of dimension $r$, 
$\{K_n(r)\ |\ n\in\Bbb Z_{\geqslant 0}\}$, give an object of the category 
$\Cal B(r)$. This object will be usually denoted by $K_{\centerdot }(r)$. 

Notice that 
two towers $K_{\centerdot }$ and 
$L_{\centerdot }$ are naturally isomorphic if $K_n=L_n$ 
for all $n\gg 0$ (all sufficiently large $n$). 
Such towers will be called almost equal.

Let $K_{\centerdot },L_{\centerdot }\in\Cal B(N)$. Then by definition 
$K_{\centerdot }\subset L_{\centerdot }$ or $L_{\centerdot }$ 
is an extension of $K_{\centerdot }$ if for all   
$m\gg 0$, $K_m\subset L_m$. $L_{\centerdot }$ is a finite 
extension of $K_{\centerdot }$ 
of degree $d=d(L_{\centerdot }/K_{\centerdot })$ if 
for all $m\gg 0$, $[L_m:K_m]=d$. 
Clearly, if $L_{\centerdot }/K_{\centerdot }$ and 
$M_{\centerdot }/L_{\centerdot }$ 
are finite extensions then $M_{\centerdot }/K_{\centerdot }$ is also finite 
and $d(M_{\centerdot }/K_{\centerdot })=d(L_{\centerdot }/K_{\centerdot })
d(M_{\centerdot }/L_{\centerdot })$.

An extension $L_{\centerdot }/K_{\centerdot }$ will be called 
separable if there is an index $m_0$ and an algebraic extension $E$ of $K_{m_0}$ 
such that $L_{\centerdot}$ is 
almost equal to 
\linebreak 
$EK_{\centerdot }:= 
\{EK_m\ |\ m\geqslant 0\}$. 
Clearly, if $L_{\centerdot }/K_{\centerdot }$ 
and $M_{\centerdot }/L_{\centerdot }$ 
are separable then $M_{\centerdot }/K_{\centerdot }$ 
is also separable. Notice also, 
that the composite of finitely many separable extensions of 
$K_{\centerdot }$ 
is again separable over $K_{\centerdot }$. Therefore, 
any finite extension $L_{\centerdot }/K_{\centerdot }$ contains 
a \lq\lq unique\rq\rq\ maximal separable over $K_{\centerdot }$ subextension 
$L_{\centerdot }^{(s)}$ (i.e. any  another maximal separable subextension  
is almost equal to $L_{\centerdot }^{(s)}$).

An extension $L_{\centerdot }/K_{\centerdot }$ will be called  purely 
inseparable if for any $n\geqslant 0$, there is an $m=m(n)\geqslant 0$ such that 
$L_n\subset K_{m}$. 
The simplest example of a purely inseparable extension of 
$K_{\centerdot }$ is 
$K'_{\centerdot }$ such that for all $m$, $K'_m=K_{m+1}$.

Suppose 
$L_{\centerdot }\supset K_{\centerdot }$ is a finite extension in $\Cal B(N)$ 
of degree $d=d(L_{\centerdot }/K_{\centerdot })$. Let 
$\widetilde{L}$ and $\widetilde{K}$ be the $p$-adic completions of the 
$\mathbin{\underset{m\geqslant 0}\to\cup}L_m$ and, resp.,  
$\mathbin{\underset{m\geqslant 0}\to\cup}K_m$. Suppose that  
$[\widetilde{L}:\widetilde{K}]=\tilde d$. Then there are the following simple properties:
\medskip 

--- $\tilde d\leqslant d$;
\medskip 

--- $\tilde d=d$ iff $L_{\centerdot }$ is separable over $K_{\centerdot }$;
\medskip 

--- $\tilde d=1$ iff $L_{\centerdot }$ is purely inseparable over $K_{\centerdot }$;
\medskip 

---  if $m_0\geqslant 0$ is such that $L_{m_0}\widetilde{K}=\widetilde{L}$ 
then $L_{\centerdot}^{(s)}=L_{m_0}K_{\centerdot }$ and $L_{\centerdot }$ is 
purely inseparable over $L_{\centerdot }^{(s)}$;
\medskip 

--- if $L_{\centerdot }^{(i)}:=\{L_m\cap\widetilde{K}\ |\ m\geqslant 0\}$ then 
$L_{\centerdot }^{(i)}$ is the maximal 
purely inseparable extension of $K_{\centerdot }$ in $L_{\centerdot }$ and $L_{\centerdot }$ 
is separable over $L_{\centerdot }^{(i)}$.
\medskip

\subsubhead {\rm 4.2.} The category   
$\Cal B^a(N)$, $N\in\Bbb Z_{\geqslant 0}$ 
\endsubsubhead

\definition{Definition} $\Cal B^a(N)$ is a full subcategory  in 
$\Cal B(N)$ consisting of  
$K_{\centerdot}\in\Cal B(N)$ such that 
there is an index $n^*=n^*(K_{\centerdot })$ and 
$c^*=c^*(n^*, K_{\centerdot })>0$ such that for all  
$n\geqslant n^*$,    
\newline a) $[K_{n+1}:K_n]=p^N$ and $\bar e(K_{n+1}/K_n)=(p,\dots ,p)\in\Bbb Z^N$; 
\newline 
b)  if $1\leqslant r\leqslant N$ then  $\pr _1j(K_{n+1}(r)/K_n(r))\geqslant p^nc^*$. 
($\pr _1(j)$ denotes the first coordinate of $j\in J_r\subset\Bbb Q^r$.)
\enddefinition 

\remark{Remark} 
\newline 
a) If $K_{\centerdot }\in\Cal B^{a}(N)$ then for 
$n\geqslant n^*(K_{\centerdot })$, all $K_n$ have the same 
last residue field. 
\newline 
b) With the above notation $K_{\centerdot }\in\Cal B^a(N)$ will 
be sometimes called a tower with parameters 
$n^*$ and $c^*$; 
notice that any ${n'}^*\geqslant n^*$ 
and $0<c'\leqslant c^*$ also can be taken as parameters for $K_{\centerdot }$. 
\endremark

\proclaim{Proposition 4.1} 
Suppose $K_{\centerdot }, L_{\centerdot }\in\Cal B(N)$ and $L_{\centerdot }$ is a finite 
separable extension of $K_{\centerdot }$. If  
$K_{\centerdot }\in\Cal B^a(N)$ then $L_{\centerdot }\in\Cal B^a(N)$. 
\endproclaim 

\demo{Proof} Suppose $K_{\centerdot }$ has parameters $n^*=n^*(K_{\centerdot })$ 
and $c^*=c^*(n^*,K_{\centerdot })$. 

If $L_{\centerdot }=\{L_m\ |\ m\geqslant 0\}$ then we 
can assume that there is an $m_0\geqslant n^*$ such that 
for all $m\geqslant m_0$, $L_{m+1}=L_{m}K_{m+1}$ and $[L_m:K_m]=d(L_{\centerdot }/K_{\centerdot })$ 
is independent on $m$. 
This implies that $[L_{m+1}:L_m]=p^N$ and $\bar e(L_{m+1}/L_m)=(p,\dots ,p)$ if $m\geqslant m_0$. 
In other words, $L_{\centerdot }$ satisfies the requirement a) of 
the above definition of objects in  $\Cal B^a(N)$.

Prove that $L_{\centerdot }$ satisfies the condition b) from 
the definition of objects from $\Cal B^a(N)$. 
By induction on $N\geqslant 0$ we can assume also that 
$L_m(N-1)=K_m(N-1)$ if $m\geqslant m_0$.

Let $\alpha _m=j(L_m/K_m)$ and $j_m=j(K_{m+1}/K_m)$.

\proclaim{Lemma 4.2} If $m\geqslant m_0$ then 
$\alpha _{m+1}\leqslant \operatorname{max}\{p\alpha _m-(p-1)j_m,\alpha _m\}$. 
\endproclaim 

\demo{Proof} By the composition property of Herbrand's function we have 
$$\varphi _{L_{m+1}/K_m}(j)=\varphi _{K_{m+1}/K_m}
\left (\varphi _{L_{m+1}/K_{m+1}}(j)\right )             \tag{1}$$
for any $j\in J(N)$. Looking at the last edge points we obtain 
$$j(L_{m+1}/K_m)=\operatorname{max}
\left\{\varphi _{K_{m+1}/K_m}(\alpha _{m+1}),j_m\right\}.$$
 On the other hand, 
$L_{m+1}=L_mK_{m+1}$ implies that 
$j(L_{m+1}/K_m)=\operatorname{max}\{\alpha _m,j_m\}$. 

Therefore, 

--- if $\alpha _m\geqslant j_m$ then $\varphi _{K_{m+1}/K_m}(\alpha _{m+1})\leqslant \alpha _m$;

--- if $\alpha _m<j_m$ then $\alpha _m$ and $\varphi _{K_{m+1}/K_m}(\alpha _{m+1})$ coincide 
because the both appear as 2nd coordinates of the prelast edge point of the $\varphi _{L_{m+1}/K_m}$.

It remains only to notice that for $j\in J_N$, 

$$\varphi _{K_{m+1}/K_m}(j)=\cases j, \text{ if } j\leqslant j_m \\ 
j_m+\frac{1}{p}(j-j_m), \text{ if } j\geqslant j_m.
\endcases $$ 
The lemma is proved.  
\enddemo 

\proclaim{Lemma 4.3} If $m\geqslant m_0$ and $\alpha _m<j_m$ then 
$\varphi _{L_m/K_m}=\varphi _{L_{m+1}/K_{m+1}}$. 
\endproclaim 

\demo{Proof} Notice first that $j(L_{m+1}/K_m)=\operatorname{max}\{\alpha _m,j_m\}=j_m$. 

Let $j_m'=j(L_{m+1}/L_m)$. Then for all $j\in J(N)$, 
$$\varphi _{L_{m+1}/K_m}(j)=\varphi _{L_m/K_m}(\varphi _{L_{m+1}/L_m}(j))\tag{2}$$
implies that $j(L_{m+1}/K_m)=\operatorname{max}\{\alpha _m, \varphi _{L_{m}/K_m}(j'_m)\}$. 
Therefore, $j_m=\varphi _{L_m/K_m}(j'_m)$. 
Notice that $\bar e(L_{m+1}/L_m)=\bar e(K_{m+1}/K_m)=(p,\dots ,p)$ 
implies that 
\linebreak 
$\varphi _{L_{m+1}/L_m}(j)=j$ and 
$\varphi _{K_{m+1}/K_m}(j)=j$ for all $\bar 0_N\leqslant j\leqslant j'_m$ and, resp., 
$\bar 0_N\leqslant j\leqslant j_m$. Therefore, 
the above relations $(1)$ and $(2)$ imply 
that for all $\bar 0_N\leqslant j\leqslant j'_m$, 
$$\varphi _{L_{m+1}/K_{m+1}}(j)=\varphi _{L_{m+1}/K_{m}}(j)=\varphi _{L_{m}/K_{m}}(j).$$
In addition, the point $(j'_m,j_m)$ is the last edge 
point of the graph of $\varphi _{L_{m+1}/K_m}$ with the minimal possible 
multiplicity $p$. 
But all edge points of $\varphi _{L_{m+1}/K_{m+1}}$ and 
$\varphi _{L_m/K_m}$ are situated in the area $j<j'_m$. 
Therefore, these functions coincide for all $j$.

The lemma is proved. 
\enddemo 

We continue the proof of our proposition. 

If $m\geqslant m_0$, then 
$\pr _1(j_m/p^m)\geqslant c^*$. Then Lemma 4.2 implies that 

$$\frac{\alpha _{m+1}}{p^{m+1}}\leqslant\operatorname{max}
\left\{ \frac{\alpha _m}{p^m}-\left (1-\frac{1}{p}\right )
(c^*,0,\dots ,0), \frac{\alpha _m}{p^{m+1}}\right\}.$$
Therefore, $\alpha _m/p^m$ tends to 0 and   
taking (if necessary) a bigger $m_0$ we can assume that for all  
$m\geqslant m_0$,  
$\pr _1(\alpha _m/p^m)<c^*$ and, therefore,   
$\alpha _m<j_m$. 
Then by Lemma 4.3 the Herbrand functions of the  extensions $L_m/K_m$ with 
$m\geqslant m_0$ coincide. 
Denote this function by $\varphi _{L_{\centerdot }/K_{\centerdot }}$ 
and use the relation $\varphi _{L_{\centerdot }/K_{\centerdot }}(j'_m)=j_m$, where $j'_m=j(L_{m+1}/L_m)$, 
from the proof of Lemma 4.3.

Because $\varphi _{L_{\centerdot }/K_{\centerdot }}$ is a piece-wise linear function 
the condition $\pr _1j_m\geqslant p^mc^*$ implies 
the existence of $0<c^*_1=c^*(m_0, L_{\centerdot })<c^*$ such that 
$\pr _1(j_m')\geqslant p^mc^*_1$ for all $m\geqslant m_0$. 

The proposition is proved. 
\enddemo 
\medskip 

\subsubhead {\rm 4.3.} The category $\Cal B^{fa}(N)$
\endsubsubhead 

4.3.1. Suppose $K_{\centerdot}\in\Cal B^a(N)$ with parameters $n^*=n^*(K_{\centerdot })$ 
and $c^*(n^*,K_{\centerdot })$.  

\definition{Definition} If indices $u_1,\dots ,u_{N}$ 
are such that  
$n^*\leqslant u_N\leqslant u_{N-1}\leqslant\dots\leqslant u_{1}$  
then $K_{u_1\dots u_N}=K_{u_1}(1)\dots K_{u_{N-1}}(N-1)K_{u_N}$. 
We usually denote this field (with its natural $F$-structure) as $K_{\bar u}$, where 
$\bar u=(u_1,\dots ,u_N)$. 
\enddefinition  

\definition{Definition} Denote by   
$\Cal B^{fa}(N)$ the full subcategory of all $K_{\centerdot }\in\Cal B^a(N)$ 
such that for some 
index parameter $\bar u^0=\bar u^0(K_{\centerdot })$,     
$K_{\bar u^0}$ has a standard $F$-structure. 
\enddefinition 

\remark{Remark} If $\bar u^0=\bar u^0(K_{\centerdot })$ is the above index 
parameter then we always assume that $n^*(K_{\centerdot })=u^0_N$.
\endremark

\proclaim{Proposition 4.4} Suppose $L_{\centerdot }\supset K_{\centerdot }$ 
is a finite extension in $\Cal B^a(N)$.  
Then there is a finite Galois extension $\widetilde{L}_{\centerdot }$ of $K_{\centerdot}$ such 
that $\widetilde{L}_{\centerdot}\supset L_{\centerdot}$ and 
$L_{\centerdot }\in\Cal B^{fa}(N)$.
\endproclaim 

\demo{Proof} Let $n^*=n^*(L_{\centerdot })=n^*(K_{\centerdot })$. 
Choose a finite Galois extension $E$ of $K_{n^*}$ such that 
$E_{\centerdot }=EK_{\centerdot }\supset L_{\centerdot }$. 
Then $E_{\centerdot }\in\Cal B^a(N)$, cf. n.4.2,  and we can assume that $n^*=n^*(E_{\centerdot })$. 
Take a finite extension $F$ of $E_{n^*}(N-1)$ such that $(E_{n^*}F,F)$ is standard in the category 
$\LC (N)$. 

Let $F_{\centerdot }=FK_{\centerdot }(N-1)$. We 
can assume that $m^*:=n^*(F_{\centerdot })=n^*(K_{\centerdot }(N-1))\geqslant n^*$. 
By induction there is a finite Galois extension $H$ of $K(N-1)_{m^*}$ such that 
$H_{\centerdot }=HK_{\centerdot }(N-1)\supset F_{\centerdot }$ 
and $H_{\centerdot }\in\Cal B^{fa}(N-1)$. Then $(E_{n^*}H,H)\in\LC (N)$ is still standard 
and, therefore, $HE_{\centerdot }\in\Cal B^{fa}(N)$. At the same time, 
$HE_{\centerdot }$ is Galois over $K_{\centerdot }$ as a composite of Galois extensions.

The proposition is proved. 
\enddemo 
\medskip 

4.3.2. 
The following proposition (or more precisely, its applications below)  
plays a crucial role in the construction of 
an analogue of the field-of-norms functor.

\proclaim{Proposition 4.5} Suppose $E_{\centerdot }\in\Cal B^{fa}(N)$. 
Then for any $u\geqslant u_N^0(E_{\centerdot })$, there is a  
$v=v(u)\geqslant{u}$ such that $(E_{u}E_{v}(N-1), E_{v}(N-1))\in\LC (N)$ 
is standard.
\endproclaim

In nn.4.3.3-4.3.6 below we assume that this proposition is proved and 
consider its applications. We need these applications later in our 
construction of  
the field-of-norms functor.  
We also need them in dimension $<N$, when proving the above Proposition 4.5   
by induction on $N$ in n.4.4.   
\medskip 

\subsubhead {\rm 4.3.3.} Functions $m_r$, $1\leqslant r<N$ 
\endsubsubhead

\proclaim{Proposition 4.6} Suppose $K_{\centerdot }\in\Cal B^{fa}(N)$ 
with the index parameter $\bar u^0(K_{\centerdot })=(u_1^0,\dots ,u_N^0)$. 
Then for 
$1\leqslant r< N$, there are non-decreasing functions 
$m_r:\Bbb Z_{\geqslant u_{r+1}^0}\longrightarrow\Bbb Z_{\geqslant u_r^0}$ such that 
for any $u_1,\dots ,u_N$ such that  
$u_{N-1}\geqslant m_{N-1}(u_N)$, \dots , $u_1\geqslant m_{1}(u_{2})$, 
 $K_{u_1u_2\dots u_{N}}$ has a standard $F$-structure. 
\endproclaim 

\demo{Proof} Use induction on $N$. 

Then for 
$K_{\centerdot }(N-1)\in\Cal B^{fa}(N-1)$, there are functions 
$m_r:\Bbb Z_{\geqslant u_{r+1}^0}\longrightarrow\Bbb Z_{\geqslant u_r^0}$, 
where $1\leqslant r\leqslant N-2$, 
such that 
if $u_{N-1}\geqslant u_{N-1}^0$,  
$u_{N-2}\geqslant m_{N-2}(u_{N-1})$,\dots , $u_1\geqslant m_{1}(u_{2})$ then 
$K(N-1)_{u_2\dots u_N}$ has a standard $F$-structure. 

If $u\geqslant u_N^0$, take $v=v(u)\geqslant u_{N-1}^0$ from Proposition 4.5.  
Then define $m_{N-1}:\Bbb Z_{\geqslant u_N^0}\longrightarrow
\Bbb Z_{\geqslant u_{N-1}^0}$ by the relation 
$$m_{N-1}(u)=\operatorname{max}\{v(u')\ |\ u^0_N\leqslant u'\leqslant u\}.$$
Then this collection of functions $m_r$, $1\leqslant r<N$, 
satisfies the requirements of our proposition.
\enddemo 

\remark{Remark} With the above notation, suppose the indices 
$(v_1^0,\dots ,v_N^0)$ are such that 
$v_1^0\geqslant\dots \geqslant v^0_N$ and the functions 
$n_r:\Bbb Z_{\geqslant v_{r+1}^0}\longrightarrow\Bbb Z_{\geqslant v_r^0}$, 
$1\leqslant r<N$, are such that $v^0_{r+1}\geqslant u^0_{r+1}$ and 
$n_r(u)\geqslant m_r(u)$ for all $u\geqslant v_{r+1}^0$. Then the proposition holds also with 
the indices $v_1^0,\dots ,v_N^0$ and the functions $n_{N-1},\dots ,n_1$. In particular, 
we can assume (if necessary) that the functions $m_r$ from our proposition 
are strictly increasing.
\endremark 

\medskip

\subsubhead {\rm 4.3.4.} Local parameters 
\endsubsubhead 

Suppose $K_{\centerdot }\in\Cal B^{fa}(N)$ and 
for $1\leqslant r <N$, $m_r:\Bbb Z_{\geqslant u_{r+1}^0}\longrightarrow 
\Bbb Z_{\geqslant u_r^0}$ are corresponding functions from the above proposition. 
We always agree to assume in this situation that $n^*(K_{\centerdot })=u_N^0$ 
and $m_r(u_{r+1}^0)=u_r^0$ for all $1\leqslant r<N$. 

Let $1\leqslant r\leqslant N$ and let indices 
$u_1,\dots ,u_r$ be such that $u_r\geqslant u_r^0$, 
$u_{r-1}=m_r(u_r)$,\dots , $u_1=m_1(u_2)$. Let $t_{u_r}^{(r)}$ be an $r$-th 
local parameter in the field 
\linebreak 
$K_{u_1}(1)K_{u_2}(2)\dots K_{u_r}(r)$. 

\proclaim{Proposition 4.7} For any indices $u_1$,\dots , $u_N$ such that 
$u_N\geqslant u_N^0$, 
\linebreak 
$u_{N-1}\geqslant m_{N-1}(u_N)$,\dots ,
$u_1\geqslant m_1(u_2)$, the above introduced elements $t_{u_1}^{(1)},\dots ,t_{u_N}^{(N)}$ 
give a system of local parameters in the field 
$K_{\bar u}=K_{u_1}(1)\dots K_{u_{N-1}}(N-1)K_{u_N}$.
\endproclaim 

\demo{Proof} If $N=1$ there is nothing to prove. 

If $N>1$ we can assume by induction that $t_{u_1}^{(1)},\dots ,t_{u_{N-1}}^{(N-1)}$ 
is a system of local parameters in $E=K_{u_1}(1)\dots K_{u_{N-1}}(N-1)$.

Let $u'_{N-1}=m_{N-1}(u_N)$, $u'_{N-2}=m_{N-2}(u'_{N-1})$,\dots ,
$u'_1=m_1(u'_2)$. Let   
$E'=K_{u'_1}(1)\dots K_{u'_{N-1}}(N-1)$ and set $K_{\bar u'}=E'K_{u_N}$. Then 
$E'\subset E$ and $(K_{\bar u'},E')\in\LC (N)$ is standard. 
Therefore,  $(K_{\bar u'}E,E)$ is also standard, i.e. 
$t_{u_N}^{(N)}$ extends the system of local parameters 
$t_{u_1}^{(1)},\dots ,t_{u_{N-1}}^{(N-1)}$ of $E$ to a system 
of local parameters of $K_{\bar u}=K_{\bar u'}E$. 

The proposition is proved. 
\enddemo 

\subsubhead{\rm 4.3.5.} Construction of special extensions 
\endsubsubhead  

Assume that $K_{\centerdot }\in\Cal B^{fa}(N)$ is given via the 
above notation. Assume in addition that the functions $m_r$, 
$1\leqslant r<N$, are strictly increasing. 

For any $n\in\Bbb Z_{\geqslant 0}$, set $v_N^n=u^0_N+n$ and define 
the vector $\bar v^n=(v_1^n,\dots ,v_N^n)$ by the relations 
$v_{N-1}^n=m_{N-1}(v^n_N+1)$,\dots , $v_1^n=m_1(v_2^n+1)$. 
Notice that 
for any indices $w_1,\dots ,w_N$ such that 
$v^n_r\leqslant w_r\leqslant v^n_r+1$ with $1\leqslant r\leqslant N$, 
the field $K_{w_1\dots w_N}$ has a standard $F$-structure. 
Indeed, for any $1\leqslant r< N$, 
$m_r(w_{r+1})\leqslant m_r(v^n_{r+1}+1)=v^n_r\leqslant w_r$.

Set for all $n\geqslant 0$, $\bar u^{n+1}=(v^n_1+1,\dots ,v^n_N+1)$,  
$\Cal O_{\bar v^n}=\Cal O_{K_{\bar v^n}}$ and 
$\Cal O_{\bar u^n}=\Cal O_{K_{\bar u^n}}$. 

Notice that we have a natural embedding $\Cal O_{\bar v^n}\subset\Cal O_{\bar u^n}$. 
In addition, $\Cal O_{\bar u^n}\subset\Cal O_{\bar v^{n+1}}$. 
Indeed, $u_N^n=v_N^n+1=v_N^{n+1}$ and if for some 
$1\leqslant r<N$, $u^n_{r+1}\leqslant v_{r+1}^{n+1}$, then 
$$u_r^n=v_r^n+1=m_r(v_{r+1}^n+1)+1=m_r(u_{r+1}^n)+1$$
$$\leqslant m_r(v_{r+1}^{n+1})+1\leqslant m_r(v_{r+1}^{n+1}+1)=v_r^{n+1}$$

For any $u\geqslant u_N^0$, 
let $v_{K_u}$ be the canonical $N$-valuation associated with $K_u$. 
Then 
$v_{K_{\centerdot }}:=v_{K_u}/p^u$ does not depend on the choice of $u$. 
We have also 1-valuations $v^1_{K_u}:=\pr _1 v_{K_u}$ and $v^1_{K_{\centerdot }}=\pr _1 v_{K_{\centerdot }}$. 
For any 
$c>0$, set 
$$\m_{K_{\centerdot }}^1(c)=\{o\in\Cal O_{\Bbb C(N)_p}\ |\ v^1_{K_{\centerdot }}(o)\geqslant c\}.$$
For any subring $O$ in $\Cal O_{\Bbb C(N)_p}$, agree to denote by 
$O\operatorname{mod}\m_{K_{\centerdot }}^1(c)$ the image of $O$ in $\Cal O_{\Bbb C(N)_p}
\operatorname{mod}\m^1_{K_{\centerdot }}(c)$. 
Notice that for any $n\geqslant 0$, there is a natural inclusion 
$\Cal O_{\bar u^n}\operatorname{mod}\m_{K_{\centerdot }}^1(c)\subset 
\Cal O_{\bar v^{n}}\operatorname{mod}\m_{K_{\centerdot }}^1(c)$. 

\proclaim{Proposition 4.8} 
Let $c_1^*=c^*(u^0_N, K_{\centerdot })/p$. Then for all $n\geqslant 0$, 
the $p$-th power map 
induces a ring epimorphism 
$$\Cal O_{\bar u^{n+1}}\operatorname{mod}\m^1_{K_{\centerdot }}(c^*_1)
\longrightarrow\Cal O_{\bar v^n}\operatorname{mod}\m^1_{K_{\centerdot }}(c^*_1).$$
\endproclaim 

\demo{Proof} Remind that 
$$\bar u^{n+1}=(u_1^{n+1},\dots ,u_N^{n+1})=(v_1^n+1,\dots ,v_N^n+1).$$
Let $1\leqslant r\leqslant N$ and let $t_{u_r}^{(r)}$ be the $r$-th local 
parameter for $K_{\bar u^{n+1}}(r)$ from n.4.3.4. It will be sufficient 
to prove that its $p$-th power is congruent modulo 
$\m ^1_{K_{\centerdot }}(c_1^*)$ to some $r$-th 
local parameter of the fierld $K_{\bar v^n}(r)$. 
By induction we can assume that $r=N$. 

Let $E=K_{\bar v^n}$, $E'=K_{\bar v^n}(t_{u^{n+1}_N}^{(N)})\subset K_{\bar u^{n+1}}$. 
Then $[E':E]=p$ and both these fields have a standard $F$-structure. 
If $\tau\in I_{E'/E}$ and $\tau\ne\id $, then 
$$v^1_{K_{u_N^{n+1}}}\left (\tau t^{(N)}_{u_N^{n+1}}-t_{u_N^{n+1}}^{(N)}\right )\geqslant p^{v_N^{n}}c^*$$
by the definition of the parameter $c^*=c^*(n^*,K_{\centerdot })$. This implies that all conjugates to 
$t_{u^{n+1}_N}^{(N)}$ over $E$ are congruent modulo 
$\m^1_{K_{\centerdot }}(c^*/p)=\m _{K_{\centerdot }}^1(c_1^*)$. Therefore, 
$p$-th power of $t^{(N)}_{u^{n+1}_N}$ is congruent modulo 
$\m _{K_{\centerdot }}^1(c_1^*)$ to the norm $N_{E'/E}\left (t^{(N)}_{u_N^{n+1}}\right )$, which 
is 
an $N$-th 
local parameter in $K_{\bar v^n}$. 

The proposition is proved. 
\enddemo 

\proclaim{Corollary 4.9} With the above notation and assumptions there is 
a field tower 
$$K_{\bar v^0}\subset K _{\bar u^1}\subset K_{\bar v^1} \subset \dots \subset 
K_{\bar v^n}\subset K_{\bar u^{n+1}}\subset\dots $$
such that all extensions $K_{\bar u^{n+1}}/K_{\bar v^n}$ satisfy the condition C from n.2.6.  
\endproclaim 

\subsubhead{\rm 4.3.6.} Modified system of local parameters 
\endsubsubhead 

As earlier, $K_{\centerdot }\in\Cal B^{fa}(N)$ together with 
the corresponding strictly increasing functions 
$m_r:\Bbb Z_{\geqslant u_{r+1}^0}\longrightarrow\Bbb Z_{\geqslant u_r^0}$ 
for $1\leqslant r<N$. 

Define 
$U(m_1,\dots ,m_{N-1})\subset\Bbb Z^N$ as the set of $\bar u=(u_1,\dots ,u_N)$ such that 
$u_N\geqslant u_N^0$, 
$u_{N-1}\geqslant m_{N-1}(u_N+1)$,\dots ,$u_1\geqslant m_1(u_2+1)$. 

\proclaim{Proposition 4.10} For all $1\leqslant r\leqslant N$ and $u\geqslant u_r^0$, 
there are $\tau ^{(r)}_u\in K_u(r)$ such that 
\newline 
{\rm a)} ${\tau _{u+1}^{(r)}}^p\equiv\tau ^{(r)}_u\operatorname{mod}m^1_{K_{\centerdot }}(c^*_1)$;
\newline 
{\rm b)} if $\bar u=(u_1,\dots ,u_N)\in U(m_1,\dots ,m_{N-1})$ then 
$\tau _{u_1}^{(1)},\dots ,\tau _{u_r}^{(r)}$ is a system of local 
parameters in $K_{\bar u}(r)$.
\endproclaim 

\demo{Proof} Use induction on $r$. Then it will be sufficient 
to define $\tau _u^{(N)}$ with $u\geqslant u_N^0$. 

Set $\tau _{u^0_N}^{(N)}=t_{u^0_N}^{(N)}$, cf. n.4.3.4. 

Then use induction on $n\geqslant 1$. Take $\tau ^{(N)}_{u_N^0+n}\in\Cal O_{\bar u^n}$ such that 
$${\tau _{u_N^0+n}^{(N)}}^p\equiv \tau ^{(N)}_{u^0_{N}+n-1}
\operatorname{mod}m^1_{K_{\centerdot }}(c^*_1)$$
Clearly,  this is an $N$-th local parameter in $K_{\bar u}$ and, therefore, it completes 
the system $\tau _{u_{N-1}}^{(N-1)},\dots ,\tau _{u_1}^{(1)}$ to a system of 
local parameters in $K_{\bar u}$.
\enddemo

\subsubhead {\rm 4.4.}  Proof of proposition 4.5 
\endsubsubhead 

Notice that there is nothing to prove 
if $N=1$ and use induction on $N$ by assuming that 
the proposition holds in dimensions $<N$. 

Therefore, we can use the result of Corollary 4.9 in dimensions $<N$. 
It remains to notice that if $K_{\bar v^n}$ is $F$-standard then $K_{v^n_N+1}K_{\bar v^n}$ is 
infernal over $K_{\bar v^n}$. So, Proposition 4.5. follows 
from the case b) of the procedure of elimination of wild ramification from n.2.3.  

\medskip

\subhead 5. Family of fields $X(K_{\centerdot })$, $K_{\centerdot }\in\Cal B^{fa}(N)$ 
\endsubhead 

\subsubhead {\rm 5.1.} Fontaine's field $R_0(N)$ 
\endsubsubhead 
 
Recall that objects $K\in\LF _0(N)$ are realised as subfields in $\Bbb C(N)_p$. 
They are closed subfields with induced 
$F$-structure and $P$-topology. 
Any $K\in\LF _0(N)$ has a canonical valuation $v_K$ of rank $N$.

Notice that if $K'\in\LF _0(N)$ then $v_{K'}=\bar\alpha v_K$ with some $\bar\alpha\in\Bbb Q^N$, 
$\bar\alpha >\bar 0$, 
and therefore, all such valuations belong to the same class of equivalent valuations. 
If $K\in\LF _0(N)$ and $v_K$ is the extension of its canonical valuation of rank $N$ 
to $\Bbb C(N)_p$ then 
$$\Cal O_{\Bbb C(N)_p}=\{o\in\Bbb C(N)_p\ |\ v_K(o)\geqslant\bar 0_N\}.$$

Set $R(N)=\mathbin{\underset{n}\to\varprojlim}(\Cal O_{\Bbb C(N)_p}\operatorname{mod}p)_n$ where connecting 
morphisms are induced by the $p$-th power map. Then $R(N)$ is an integral domain and its fraction field 
$R_0(N)$ is a perfect field of characteristic $p$.  
The $F$-structure on $\Bbb C(N)_p$ induces an $F$-structure 
on 
$R_0(N)$ given by the decreasing sequence of subfields 
$$R_0(N)\supset R_0(N-1)\supset \dots \supset R_0(1)\supset R_0(0).$$
In addition, the field $R_0(0)$ consists of the sequences 
$\{\alpha ^{p^{-n}}\}_{n\geqslant 0}$, where $\alpha\in\bar\Bbb F_p$. The map 
$\{\alpha ^{p^{-n}}\}_{n\geqslant 0}\mapsto\alpha $ identifies $R_0(0)$ with $\bar\Bbb F_p$, 
in particular, any finite field of characteristic $p$ can be embedded naturally into $R_0(N)$. 

Notice that $R=R(1)$ and $\operatorname{Frac}R=R_0(1)$ are original notations introduced for the 
corresponding 1-dimensional objects by J.-M. Fontaine. 

Let $K_{\centerdot }\in\Cal B^a(N)$. It determines a valuation of rank $N$ on $\Bbb C(N)_p$ given by 
the formula $v_{K_{\centerdot }}=\mathbin{\underset{n\to\infty }\to\lim}(v_{K_n}/p^n)$. 
This determines the valuation $v_{R,K_{\centerdot }}$ of rank $N$ on $R_0(N)$ such that 
if $\bar r=(r_n)_{n\geqslant 0}\in R(N)$ then 
$$v_{R,K_{\centerdot }}(\bar r)=\lim _{n\to\infty }p^nv_{K_{\centerdot }}(\hat r_n)=
\lim_{n\to\infty }v_{K_n}(\hat r_n)$$
where $\hat r_n\in\Cal O_{\Bbb C(N)_p}$ is such that $\hat r_n\operatorname{mod}p=r_n$. 

Notice that if $L_{\centerdot }\in\Cal B^a(N)$ then $v_{R,L_{\centerdot }}=\bar\alpha v_{R,K_{\centerdot }}$ 
with $\bar\alpha\in\Bbb Q^N$, $\bar\alpha >\bar 0$. Therefore, 
the equivalence class of such valuations does not depend on the choice of $K_{\centerdot }$. 

If $c>0$, then (as earlier) 
$$\m^1_{R,K_{\centerdot }}(c)=\{o\in R(N)\ |\ v^1_{R,K_{\centerdot }}(o)\geqslant c\}$$ 
where $v^1_{R,K_{\centerdot }}=\pr _1v_{R,K_{\centerdot }}$. 

The following proposition is just an easy consequence of the above definitions. 

\proclaim{Proposition 5.1} For any $c>0$ such that $p\in\m _{K_{\centerdot }}^1(c)$, 
\newline 
{\rm a)} $R(N)=\varprojlim (\Cal O_{\Bbb C(N)_p}\operatorname{mod}\m^1_{K_{\centerdot }}(c))$, where 
connecting morphisms are induced by the $p$-th power map;
\newline 
{\rm b)} for any $u\geqslant 0$, the $u$-th projection 
$\pr _u:R(N)\longrightarrow \Cal O_{\Bbb C(N)_p}\operatorname{mod}\m^1_{K_{\centerdot }}(c)$ 
induces a ring identification of $R(N)\operatorname{mod}\m^1_{R,K_{\centerdot }}(p^uc)$ and 
$\Cal O_{\Bbb C(N)_p}\operatorname{mod}\m^1_{K_{\centerdot }}(c)$.
\endproclaim 

\remark{Remark} $\Cal O_{\Bbb C(N)_p}$ is equipped with the $P$-topology 
given by the inductive limit 
of $P$-topologies on all $K\in\LF _0(N)$. 
This topology induces the $P$-topology on $R(N)$ and $R_0(N)$.
\endremark 
\medskip 

\subsubhead {\rm 5.2.} The family of fields $X(K_{\centerdot})$ 
\endsubsubhead 

Suppose $K_{\centerdot }\in\Cal B^{fa}(N)$ with the parameters  
$\bar u^0(K_{\centerdot })=(u_1^0,\dots ,u_N^0)$ and $c^*=c^*(u^0_N, K_{\centerdot })$. 
As earlier in n.4.3.3, choose for all $1\leqslant r<N$ the corresponding strictly increasing functions 
$m_r:\Bbb Z_{\geqslant u^0_{r+1}}\longleftarrow\Bbb Z_{\geqslant u^0_r}$ 
and $r$-th local parameters $\tau _u^{(r)}\in K_u(r)$, where 
$u\geqslant u_r^0$. 

Set (as earlier, $c_1^*=c^*/p$)

$$\tau ^{(r)}=(\tau _u^{(r)}\operatorname{mod}\m^1_{K_{\centerdot }}(c^*_1))_{u\geqslant u_r^0}
\in\varprojlim (\Cal O_{\Bbb C(N)_p}\operatorname{mod}\m ^1_{K_{\centerdot }}(c^*_1))=R(N).$$ 

Let $k=k(K_{\centerdot })$ be the last residue field of $K_{u_N^0}$ (this is also the residue 
field for all $K_u$ with $u\geqslant u_N^0$). As it was mentioned in n.5.1, $k$ can be naturally 
indentified with a subfield in $R_0(0)\subset R_0(N)$. 

\proclaim{Proposition 5.3}
The correspondences 
$T_1\mapsto\tau ^{(1)}$,\dots , $T_N\mapsto\tau ^{(N)}$ determine a unique 
continuous embedding of the $N$-dimensional local field $k((T_N))\dots ((T_1))$ into $R_0(N)$. 
Its image is an $N$-dimensional local subfield $\Cal K$ in $R_0(N)$ with the system of local 
parameters $\tau ^{(1)},\dots ,\tau ^{(N)}$. 
\endproclaim  

\demo{Proof} We need the following obvious lemma.

\proclaim{Lemma 5.4} Suppose $L\in\LF _0(N)$ has a standard $F$-structure, 
which is compatible with given local 
parameters $t_1,\dots ,t_N$. Let $c>0$ and 
$\m^1_L(c)=\{o\in\Cal O_{\Bbb C(N)_p}\ |\ \pr _1v_L(o)\geqslant c\}$. 
Then any $o\in\Cal O_L$ can be uniquely presented modulo $\m^1_L(c)$  
in the form 
$$\sum\Sb a_1<c\endSb [\alpha _{\bar a}]t_1^{a_1}\dots t_N^{a_N}.$$
\endproclaim 
\medskip 
\remark{Remark} The coefficients $[\alpha _{\bar a}]$ are the Teichmuller representatives of 
the elements of the last residue field of $L$ and satisfy the standard restrictions from the 
beginning of n.1.1. 
\endremark 
\medskip 

Continue the proof of proposition 5.3. 

We prove first that the power series 
$$\sum\Sb \bar a\geqslant\bar 0_N\endSb \alpha _{\bar a}\tau ^{(1)a_1}\dots \tau ^{(N)a_N}\tag{3}$$
converges in $R(N)$ 
if its coefficients $\alpha _{\bar a}$ satisfy the restrictions described in n.1.1.  This is equivalent to the fact that 
for all $u\geqslant u^0_N$, the series  
 
$$\sum\Sb \bar a\geqslant\bar 0_N\endSb [\alpha _{\bar a}]^{p^{-u}}
\tau _u^{(1)a_1}\dots \tau _u^{(N)a_N}\tag{4}$$
converge to elements  $f_u\in\Cal O_{\Bbb C(N)_p}$ such that 
$f_u^p\equiv f_{u+1}\operatorname{mod}\m_{K_{\centerdot }}^1(c^*_1)$.

Let $\bar u^n=(u_1^n,\dots u_{N-1}^n,u_N^n)$ with $u_N^n=u$. Then 
for $1\leqslant r\leqslant N$, it holds $u\leqslant u_r^n$ and 
$$\tau _u^{(r)}\equiv \tau _{u_r^n}^{(r)p^{u_r^n-u}}\operatorname{mod}\m_{K_{\centerdot }}^1(c^*_1)$$ 
This means that the above series $(4)$ can be expressed in terms of 
local parameters of the field $K_{\bar u^n}$, its coefficients 
$[\alpha _{\bar a}]^{p^{-u}}$ satisfy the restrictions from 
n.1.1 and, therefore, these series converge in 
$\Cal O_{\bar u^n}\subset\Cal O_{\Bbb C(N)_p}$. 

Then  the uniqueness property from lemma 5.4 implies that 
 $$f_{u+1}^p\equiv f_u\operatorname{mod}\m _{K_{\centerdot }}^1(c^*_1)$$  
and the series $(3)$ converges in $R(N)$. 

Even more, Lemma 5.4 implies that any element from $R(N)$ can be presented 
in at most one way as a 
sum of the series $(3)$. So, the image $\Cal K$ of $k((T_N))\dots ((T_1))$ 
is an $N$-dimensional local field with the set of local parameters 
$\tau ^{(1)},\dots ,\tau ^{(N)}$. 

The proposition is proved. 
\enddemo 

Notice that the above fields $\Cal K\subset R_0(N)$ are not uniquely determined by a given 
$K_{\centerdot }\in\Cal B^{fa}(N)$. They depend also on the choice of 
functions $m_1,\dots ,m_{N-1}$ and the choice of compatible systems of local parameters 
$\{\tau _u^{(r)}\}_{u\geqslant u_r^0}$, $1\leqslant r\leqslant N$. 
Denote by $X(K_{\centerdot }; m_1,\dots ,m_{N-1})$ the family of 
all subfields $\Cal K$ which can be constructed for a given tower 
by the use of given invariants $\bar u^0(K_{\centerdot })$ together 
with an appropriate choice of strictly increasing functions $m_1,\dots ,m_{N-1}$. 
Notice that taking a bigger invariant $\bar u^0(K_{\centerdot })$ together with the 
contraction of the domain  of definition of functions $m_1,\dots ,m_{N-1}$ doesn't affect 
this family. Clearly, for a given tower $K_{\centerdot }$, the sets 
$X(K_{\centerdot }; m_1,\dots ,m_{N-1})$ form an inductive system. Its 
inductive limit 
will be denoted by $X(K_{\centerdot })$.

\subsubhead{\rm 5.3.} The categories $\LF _R(N)$ and $\widetilde{\LF}_R(N)$
\endsubsubhead 

Consider the category $\LF _R(N)$ of all $N$-dimensional closed subfields $\Cal K$ in $R_0(N)$ 
together with the induced $F$-structure given by the subfields of   
$r$-dimensional constants $\Cal K(r)=R_0(r)\cap\Cal K$, $0\leqslant r\leqslant N$. 
If $\Cal K,\Cal L\in\LF _R(N)$ then $\Hom _{\LF _R(N)}(\Cal K,\Cal L)$ consists of 
compatible with $F$-structure and $P$-continuous morphisms $f:R_0(N)\longrightarrow R_0(N)$ 
such that $f(\Cal K)\subset \Cal L$. 

Suppose that $v^1$ is a 1-dimensional valuation coinciding with one of 
(equivalent valuations) $\pr _1v_{R,K_{\centerdot }}$, where 
$K_{\centerdot}\in\Cal B^{fa}(N)$. 

For any $v^1$-adic closed subfield $\Cal L$ in $R_0(N)$ denote by $\Cal R(\Cal L)$ the $v^1$-adic closure 
of the maximal inseparable extension of $\Cal L$ in $R_0(N)$.

\definition{Definition} If $\Cal K,\Cal L\in\LF _R(N)$ then $\Cal K\sim\Cal L$ if 
for $1\leqslant r\leqslant N$, $\Cal K(r)\Cal R(\Cal K(r-1))=\Cal L(r)\Cal R(\Cal L(r-1))$, where  
the composite is taken in the category of $v^1$-adic closed subfields of $R_0(N)$.
\enddefinition 

Clearly, the above defined relation $\sim $ is an equivalence relation. 
Denote by $\widetilde{\LF }_R(N)$ the category such that its 
objects are equivalence 
classes $\cl (\Cal K)$ of all $\Cal K\in\LF _R(N)$ and 
for any $\cl (\Cal K),\cl (\Cal L)\in\widetilde{\LF }_R(N)$, 
$\Hom_{\widetilde{\LF }_R(N)}(\cl (\Cal K),\cl (\Cal L))$ consists of 
compatible with $F$-structure and $P$-continuous 
field morphisms $f:R_0(N)\longrightarrow R_0(N)$ 
such that for any $1\leqslant r\leqslant N$, 
$f(\Cal K(r))\subset \Cal L(r)\Cal R(\Cal L(r-1))$.

\remark{Remark} The usual \lq\lq 1-dimensional\rq\rq\  Krasner's Lemma implies that:
\newline 
{\it --- if $\Cal L_1,\Cal L\in\LF _R(N)$,  
$[\Cal L_1:\Cal L]=m$ and $\Cal L'\sim \Cal L$  then there is a unique $\Cal L_1'\in\LF _R(N)$ 
such that 
$\Cal L_1'\sim \Cal L_1$ and $\Cal L_1'$ is an extension of $\Cal L'$ of degree $m$.}
\newline 
In particular, we can use the concepts of finite algebraic, 
separable, Galois and purely inseparable extensions in 
$\widetilde{\LF}_R(N)$. 
\endremark 
\medskip 

\subsubhead{\rm 5.4.} Identification of elements from $X(K_{\centerdot })$, 
$K_{\centerdot }\in\Cal B^{fa}(N)$
\endsubsubhead 

\proclaim{Proposition 5.5} Suppose $K_{\centerdot }\in\Cal B^{fa}(N)$. 
Then any two elements from $X(K_{\centerdot })$ represent the same object
in $\widetilde{\LF }_R(N)$. 
\endproclaim

\demo{Proof} Let $\bar u^0(K_{\centerdot })=(u_1^0,\dots ,u_N^0)$ and 
$c^*_1=c^*(u^0_N, K_{\centerdot })/p$. 

Suppose $\Cal K\in X(K_{\centerdot })$ is obtained via a choice of 
strictly increasing functions 
$m_r:\Bbb Z_{\geqslant u_{r+1}^0}\longrightarrow\Bbb Z_{\geqslant u_r^0}$,  
and a special system of local parameters 
$\tau _u^{(r)}$, $1\leqslant r<N$, $u\geqslant u^0_r$. 

Take some $u\geqslant u_N^0$ and choose 
$\bar u=(u_1,\dots ,u_{N-1},u)\in U(m_1,\dots ,m_{N-1})$. 

Set $\Cal K_{\bar u}=\Cal K\left (\sigma ^{u-u_{N-1}}\Cal K(N-1)\right )\dots 
\left (\sigma ^{u-u_1}\Cal K(1)\right )$. (Here $\sigma $ is as usually the $p$-th power map.)
Then 
$$\sigma ^{u-u_1}\tau ^{(1)}, \dots ,\sigma ^{u-u_{N-1}}\tau ^{(N-1)}, \tau ^{(N)}$$ 
is a system of local parameters in $\Cal K_{\bar u}$ which is compatible with a given (standard) 
$F$-structure of $\Cal K_{\bar u}$. 

It is easy to see that for $1\leqslant r\leqslant N$, 
the correspondences $\sigma ^{u-u_r}\tau ^{(r)}\mapsto\tau _{u_r}^{(r)}$ 
give the identification 
$$\psi _{\bar u}:\Cal O_{\Cal K_{\bar u}}\operatorname{mod}\m^1_{R,K_{\centerdot }}(p^uc^*_1)
\simeq \Cal O_{\bar u}\operatorname{mod}\m^1_{K_{\centerdot }}(c^*_1).$$

If $\bar u'=(u_1',\dots ,u_{N-1}',u)\in U(m_1,\dots ,m_{N-1})$ 
is such that $u'_r\geqslant u_r$ for all $1\leqslant r\leqslant N$, 
then $\psi _{\bar u}$ and $\psi _{\bar u'}$ are compatible via 
natural inclusions $\Cal K_{\bar u}\subset\Cal K_{\bar u'}$ and 
$\Cal O_{\bar u}\subset\Cal O_{\bar u'}$. Therefore, the $u$-th projection 
$\pr _u:R(N)\longrightarrow\Cal O_{\Bbb C(N)_p}\operatorname{mod}m^1_{K_{\centerdot }}(c^*_1)$ 
induces the identification 
$$\psi _u:\Cal O_{\Cal K\Cal R(\Cal K(N-1))}\longrightarrow \Cal O^{(u)}
\operatorname{mod}m^1_{K_{\centerdot }}(c^*_1)$$
where $\Cal O^{(u)}$ is the valuation ring of the composite of all 
$K_{\bar u}$ with $\bar u$ running over the set of all $\bar u=(u_1,\dots ,u_{N-1},u_N)$ 
such that $u_N=u$. 

In order to understand the relation between different $\psi _u$, 
notice that if 
$\bar u'=(u_1',\dots ,u_{N-1}',u+1)\in U(m_1,\dots ,m_{N-1})$,  
then $\bar u=(u_1'-1,\dots ,u_{N-1}'-1,u)\in U(m_1,\dots ,m_{N-1})$ 
(because the functions $m_r$ are strictly increasing) and 
$\Cal K_{\bar u}=\Cal K_{\bar u'}$. This implies that $\psi _{\bar u}$ and 
$\psi _{\bar u'}$  fit into a commutative diagram via 
the natural projection 
$$\Cal O_{\Cal K_{\bar u'}}\operatorname{mod}\m^1_{R,K_{\centerdot }}(p^{u+1}c^*_1)
\longrightarrow 
\Cal O_{\Cal K_{\bar u}}\operatorname{mod}\m^1_{R,K_{\centerdot }}(p^{u}c^*_1)$$
and the restriction of the transition morphism of the projective system 
\linebreak 
$\Cal O_{\Bbb C(N)_p}\operatorname{mod}\m^1_{K_{\centerdot }}(c^*_1)$ from the definition of $R(N)$. 
Therefore, $\varprojlim\psi _u$ identifies  
\linebreak 
$\Cal O_{\Cal K\Cal R(\Cal K(N-1))}$ with  
$\varprojlim \Cal O^{(u)}\operatorname{mod}\m^1_{K_{\centerdot }}(c^*_1)\subset R(N)$. 
In particular, $\Cal K\Cal R(\Cal K(N-1))$ does not depend on the choice of 
$\Cal K$ in $X(K_{\centerdot })$. 

The proposition is proved. 
\enddemo 
\medskip 

5.5. Let $K_{\centerdot },L_{\centerdot }\in\Cal B^{fa }(N)$, 
$\Cal K\in X(K_{\centerdot})$ and $\Cal L=X(L_{\centerdot })$. 
Let $\widetilde K$ and $\widetilde{L}$ be the $p$-adic completions of 
$\mathbin{\underset{m\geqslant 0}\to\cup}K_m$ and, resp.,  
$\mathbin{\underset{m\geqslant 0}\to\cup}L_m$. 
Notice that if $\bar\Cal K$ is purely inseparable over $\bar\Cal L$ then 
$\widetilde{K}=\widetilde{L}$. Inversely, we have the following property.

\proclaim{Proposition 5.6} 
With the above notation, if $\widetilde{K}=\widetilde{L}$ 
then $\Cal R(\Cal K)=\Cal R(\Cal L)$. 
\endproclaim

The proof is straightforward.  

\medskip 
\medskip

\subhead 6. Separable extensions in $\Cal B^{fa}(N)$ and $\widetilde{\LF }_R(N)$ 
\endsubhead 
\medskip 

6.1. In this subsection we prove that the correspondence $K_{\centerdot }\mapsto\cl (\Cal K)$, 
where $\Cal K\in X(K_{\centerdot })$,  
transforms finite separable extensions in $\Cal B^{fa}(N)$ to finite 
separable extensions of the same degree in $\widetilde{\LF }_R(N)$. 

\proclaim{Proposition 6.1} Suppose $L_{\centerdot }, K_{\centerdot }\in\Cal B^{fa}(N)$ and  
$L_{\centerdot }\supset K_{\centerdot }$ is separable and finite  
of degree $d(L_{\centerdot }/K_{\centerdot })=d$. 
Then for any $\Cal K\in X(K_{\centerdot })$, there is an 
$\Cal L\in X(L_{\centerdot })$ such that $\Cal L$ is a separable extension of 
$\Cal K$ of degree $d$.
\endproclaim 

\demo{Proof} By induction on $N\geqslant 0$, we can assume that 
$K_{\centerdot }(N-1)=L_{\centerdot }(N-1)$.

We can assume also that: 
\medskip 

--- $\bar u^0(K_{\centerdot })=\bar u^0(L_{\centerdot })=(u^0_1,\dots ,u_N^0)$;
\medskip 

--- $c^*(u^0_N, K_{\centerdot })=c^*(u^0_N, L_{\centerdot })=pc^*_1$;
\medskip

--- there are strictly increasing functions 
$m_r:\Bbb Z_{\geqslant u^0_{r+1}}\longrightarrow\Bbb Z_{\geqslant u^0_r}$, 
where $1\leqslant r<N$,  such that $m_r(u^0_{r+1})=u_r^0$ and 
if $\bar u\in U(m_1,\dots ,m_{N-1})$ then both $L_{\bar u}$ and $K_{\bar u}$ have a standard 
$F$-structure; 
\medskip 

--- $\Cal K\in X(K_{\centerdot }; m_1,\dots ,m_{N-1})$; 
\medskip 

--- for all $u\geqslant u^0_N$, the Herbrand functions of extensions $L_u/K_u$ 
coincide and are equal to  
$\varphi _{L_{\centerdot }/K_{\centerdot }}$; 
\medskip

---  the initial choice of $u_N^0$ provides us with the inequality $\pr _1(j(L_{\centerdot }/K_{\centerdot }))
+\delta _{1N}<p^{u^0_N}c^*_1/2$ (here and everywhere below $\delta _{1N}$ is the Kronecker symbol). 
\medskip 

As usually, we denote by 
$(i(L_{\centerdot }/K_{\centerdot }),j(L_{\centerdot }/K_{\centerdot }))$ the last edge point of the graph of 
$\varphi _{L_{\centerdot }/K_{\centerdot }}$  
and use the notation $\pr _1(i(L_{\centerdot }/K_{\centerdot }))=i^1$,  
$\pr _1(j(L_{\centerdot }/K_{\centerdot }))=j^1$, 
$\pr _1(\bar e(L_{\centerdot }/K_{\centerdot }))=e^1$.   

Consider the corresponding sequence of multi-indices 
$\bar u^0,\bar v^0,\dots ,\bar u^n,\bar v^n,\dots $ from n.4.3.5 and the corresponding 
field towers

$$L_{\bar v^0}\subset L_{\bar u^1}\subset L_{\bar v^1}\subset\dots 
\subset L_{\bar u^n}\subset L_{\bar v^n}\subset \dots $$

$$K_{\bar v^0}\subset K_{\bar u^1}\subset K_{\bar v^1}\subset\dots 
\subset K_{\bar u^n}\subset K_{\bar v^n}\subset \dots $$

For any $u\geqslant u_N^0$, set $n=n(u)=u-u^0_N$. So, 
$\bar u^n=(u_1^n,\dots ,u_{N-1}^n,u)$.

Consider $a_{iu}\in \Cal O_{K_{\bar u^n}}$, 
where $1\leqslant i\leqslant d$ 
and $u\geqslant u^0_N$, such that 
\medskip 

--- there is $N$-th local parameter $\eta _{u^0_N}$ in $L_{\bar u^0}$ 
such that 
$$\eta _{u^0_N}^d+a_{1u^0_N}\eta _{u^0_N}^{d-1}+\dots +a_{du^0_N}=0;$$ 
\medskip 

--- for all $1\leqslant i\leqslant d$ and $u\geqslant u^0_N$, 
$a_{iu}\equiv a_{i,u+1}^p\operatorname{mod}\m^1_{K_{\centerdot }}(c^*_1)$ or, 
equivalently, 
there are $\alpha _i\in\Cal O_{\Cal K}$ such that 
$\pr _u(\alpha _i)=a_{iu}\operatorname{mod}m^1_{K_{\centerdot }}(c^*_1)$, 
where $\pr _u$ is the projection from  
$R(N)=\mathbin{\underset{u}\to\varprojlim }(\Cal O_{\Bbb C(N)_p}
\operatorname{mod}m^1_{K_{\centerdot }}(c^*_1))_u$ 
to its $u$-th component. 
\medskip

Let $F_u(T)=T^d+a_{1u}T^{d-1}+\dots +a_{du}$, where $u\geqslant u^0_N$. Then 
all $F_u(T)$ are $N$-th Eisenstein polynomials in 
$\Cal O_{K_{\bar u^n}}[T]$ 
(i.e. their images in $\Cal O_{K^{(N-1)}_{\bar u^n}}[T]$ are Eisenstein polynomials, 
where $K_{\bar u^n}^{(N-1)}$ is the pre-last residue field of $K_{\bar u^n}$) 
and $F_{u^0_N}(\eta _{u^0_N})=0$. 
We want to prove that for 
$c^*_2=e^1c^*_1/2$, there are $N$-th local parameters 
$\eta _u\in L_{\bar u^n}$ with $u>u^0_N$, 
such that for all $u\geqslant u^0_N$, 
$$\eta _{u+1}^p-\eta _u\in\m^1_{L_{u}}(p^uc^*_2)$$
(notice that $\eta _{u^0_N}$ has been chosen earlier).

Suppose $u\geqslant u^0_N$ and we have already constructed such 
elements $\eta _v$ for all $v$ such that $u^0_N\leqslant v\leqslant u$. 

\proclaim{Lemma 6.2} If $\theta _{u+1}\in\Cal O_{\Bbb C(N)_p}$ is a root of $F_{u+1}(T)$ then there is 
a unique root $\theta _u\in\Cal O_{\Bbb C(N)_p}$ of $F_u(T)$ such that 
$\theta _u-\theta _{u+1}^p\in m^1_{L_u}(p^uc^*_2)$. 
\endproclaim

\demo{Proof of lemma} 
 Clearly, $F_u(\theta _{u+1}^p)\in m^1_{K_u}(p^uc^*_1)$. 
Let $v_{K_u}(F_u(\theta ^p))=j_u+(0,\dots ,0,1)$. Then by assumptions from n.6.1 
$$\pr _1(j_u)\geqslant p^uc^*_1-\delta _{1N}> 2j^1+\delta _{1N}\geqslant j^1.$$
Therefore, $j_u>j(L_u/K_u)=j(L_{\centerdot }/K_{\centerdot })$ and we can apply Krasner's lemma, cf. n.3.4.   
This lemma gives the existence of a unique root $\theta _u\in\Cal O_{\Bbb C(N)_p}$ of $F_u(T)$ such that 
$v_{L_u}(\theta _{u+1}^p-\theta _u)=i_u+(0,\dots ,0,1)$ with 
$\varphi _{L_{\centerdot }/K_{\centerdot }}(i_u)=j_u$. Because $j_u\geqslant j(L_u/K_u)$, 
we have 
$$\frac{j_u-j(L_u/K_u)}{i_u-i(L_u/K_u)}=\bar e^{-1}(L_u/K_u)$$
and this implies  
$$\pr _1(i_u)\geqslant e^1(\pr _1j_u-j^1)\geqslant e^1(p^uc^*_1-\delta _{1N}-j^1)
>e^1(p^uc^*_1-\frac{1}{2}p^{u^0_N}c^*_1)\geqslant\frac{1}{2}e^1p^uc^*_1=p^uc^*_2.$$

The lemma is proved.
\enddemo 

Notice that $L_{\bar u^{n+1}}=L_{\bar u^n}K_{\bar u^{n+1}}=K_{\bar u^{n+1}}(\theta _u)$ is of degree 
$d$ over $K_{\bar u^{n+1}}$ and, therefore, $F_u(T)$ is still irreducible 
over $K_{\bar u^{n+1}}$. Therefore, there is a 
\linebreak 
$\tau \in\Gamma _{K_{\bar u^{n+1}}}$ 
such that $\tau (\theta _u)=\eta _u$. 
Take $\eta _{u+1}=\tau (\theta _{u+1})$. 
Then the uniqueness of $\theta _u$ in the above lemma  
implies that the field $K_{\bar u^{n+1}}(\eta _{u+1})$ contains the field 
$K_{\bar u^{n+1}}(\eta _u)=L_{\bar u^{n+1}}$.  Therefore, these fields are equal, because 
they are both of the same degree $d$ over $K_{\bar u^{n+1}}$. 

Finally, $\eta _{u+1}$ is an $N$-th local parameter in 
$L_{\bar u^{n+1}}$ because $F_{u+1}$ is an $N$-th Eisenstein polynomial in 
$\Cal O_{K^{n+1}}[X]$ and the existence of the sequence $\eta _u$, $u\geqslant u^0_N$, 
is proved.

Let 
$$\eta =\mathbin{\underset{u}\to\varprojlim } \eta _u 
\in\mathbin{\underset{u}\to\varprojlim}\left (\Cal O_{\Bbb C(N)_p}
\operatorname{mod}m^1_{L_{\centerdot}}(c^*_2)\right )_u=R(N)$$
Then $\eta $ is a root of $N$-th Eisenstein polynomial 
$$F(T)=T^d+\alpha _1T^{d-1}+\dots +\alpha _d\in\Cal O_{\Cal K}[T].$$
Therefore, $\Cal L=\Cal K(\eta )$ is of degree $d$ over $\Cal K$, $\Cal L$ 
has standard $F$-structure and $\eta $ is its $N$-th local parameter. Clearly, 
$\Cal L\in X(L_{\centerdot };m_1,\dots ,m_{N-1})$. 

We now prove that $\Cal L$ is separable over $\Cal K$. Indeed, notice first that 
\medskip 

a) any other root of $F_{u^0_N}$ equals $\tau\eta _{u^0_N}$ for a suitable 
automorphism $\tau $ of $\Bbb C(N)_p$ such that 
$\tau |_{K_{\bar u^n}}=\id $ for all $n\geqslant 0$; 
\medskip 

b) with the above notation, for a sufficiently large $u$, 
$i^1+\delta _{1N}<p^uc^*_2$ and, therefore, 
$\tau\eta _u\not\equiv\eta _u\operatorname{mod}\m ^1_{L_{\centerdot }}(c^*_2)$. 
\medskip 

Therefore, $\tau\eta :=\varprojlim\tau\eta _u\in R(N)$ is again a root of $F(X)$ 
which is different from $\eta $. Therefore, $F(X)$ has $d$ distinct roots in $R_0(N)$. 

The proposition is proved. 
\enddemo

\proclaim{Corollary 6.3} Under assumptions from the above proposition: 
\newline 
{\rm a)} 
there is a natural identification of the set of all isomorphic embeddings 
$\iota $ of $L_{\centerdot }$ into $\Bbb C(N)_p$
such that $\iota |_{K_{\centerdot }}=\id $ 
and the set of all isomorphic embeddings $\iota :\Cal L\longrightarrow R_0(N)$ such that 
$\iota |_{\Cal K}=\id $;
\newline 
{\rm b)} $\varphi _{L_{\centerdot }/K_{\centerdot }}=\varphi _{\Cal L/\Cal K}$.
\endproclaim  
\medskip 

6.2. With the above notation we are going to prove now that 
for a sufficiently large separable extension $E_{\centerdot}$ of $K_{\centerdot }$, 
the appropriate $\Cal E\in\Cal X(E_{\centerdot})$ contains any 
given separable extension of $\Cal K$ in $R_0(N)$.

\proclaim{Proposition 6.4} Suppose $K_{\centerdot }\in\Cal B^{fa}(N)$, 
$\Cal K\in X(K_{\centerdot })$ and $\Cal L$ is a finite separable extension of 
$\Cal K$ with standard $F$-structure such that $\Cal K(N-1)=\Cal L(N-1)$. 
Then there is an $L_{\centerdot }\in\Cal B^{fa}(N)$ 
and a field embedding $\iota :\Cal L\longrightarrow R_0(N)$ such that 
\newline 
{\rm a)} $L_{\centerdot }$ is a separable extension of $K_{\centerdot }$ 
of degree $d=[\Cal L :\Cal K]$;
\newline 
{\rm b)} $\iota (\Cal L)\in X(L_{\centerdot })$.
\endproclaim 

\demo{Proof} We can assume that:
\medskip 

--- there are parameters $\bar u^0(K_{\centerdot })=(u^0_1,\dots ,u^0_N)$, 
$c^*(u^0_N, K_{\centerdot })=pc^*_1$ and strictly increasing functions 
$m_r:\Bbb Z_{\geqslant u^0_{r+1}}\mapsto\Bbb Z_{\geqslant u^0_r}$, where $1\leqslant r<N$, 
such that 
$\Cal K\in X(K_{\centerdot }; m_1,\dots ,m_{N-1})$;
\medskip 

--- $\Cal O_{\Cal L}=\Cal O_{\Cal K}[\theta ]$ 
where $\theta $ is a root of $N$-th Eisenstein polynomial 
$\Cal F(T)=T^d+\alpha _1T^{d-1}+\dots +\alpha _d\in\Cal O_{\Cal K}[T]$;
\medskip 

--- $2v^1_{\Cal K}(D(\Cal F))<p^{u^0_N}c^*_1-1$, where $D(\Cal F)$ is the discriminant 
of $\Cal F$ over $\Cal O_{\Cal K}$. 
\medskip 

As earlier consider the sequence $\bar u^0=u^0(K_{\centerdot }), \bar u^1,\dots ,\bar u^n,\dots $ 
and set $u=n+u^0_N$. For $u\geqslant u^0_N$, introduce 
the polynomials 
$$F_u(T)=T^d+a_{1u}T^{d-1}+\dots +a_{du}\in\Cal O_{K_{\bar u^n}}[T]$$
where $a_{iu}\operatorname{mod}\m^1_{K_{\centerdot }}(c^*_1)=\pr _{\bar u^n}(\alpha _{in})$ 
for $1\leqslant i\leqslant d$. Notice that for $u\geqslant u^0_N$, 
$$D(F_u)\operatorname{mod}\m ^1_{K_{\centerdot }}(c^*_1)=\pr _{\bar u^n}D(\Cal F)\ne 0.$$ 

For $u\geqslant u^0_N$, we will prove the existence of roots 
$\eta _u\in \Cal O_{\Bbb C(N)_p}$ of $F_u(T)$ 
such that if $M_u=K_{\bar u^n}(\eta _u)$, then 
\medskip 
 
1) $\eta _{u^0_N}$ is $N$-th local parameter in $M_{u^0_N}$;
\medskip 

2) $\eta _u-\eta _{u+1}^p\in\m^1_{M_u}(p^uc^*_1/2)$.  
\medskip 

Suppose such roots $\eta _{u^0_N},\dots ,\eta _u$ have been already constructed. 

Let $\pr _1\bar e(K_{\bar u^n}/K_u)=e_{1u}$.

Let $\theta _{u+1}\in \Cal O_{\Bbb C(N)_p}$ be a root of $F_{u+1}(T)$. Then 
$F_u(\theta _{u+1}^p)\in \m ^1_{K_{\bar u^n}}(e_{1u}p^uc^*_1)$ and, 
therefore,   
$v_{K_{\bar u^n}}(F_u(\theta _{u+1}^p))=j_u+(0,\dots ,0,1)$ with 
$$\pr _1(j_u)+\delta _{1N}\geqslant e_{1u}p^uc^*_1.$$

\proclaim{Lemma 6.5} $j_u>j(M_u/K_{\bar u^n})$. 
\endproclaim 

\demo{Proof of lemma} From Corollary 3.5 we have   

$$j(M_u/K_{\bar u^n})\leqslant 2v_{K_{\bar u^n}}(D(F_u))=
2\bar e(K_{\bar u^n}/K_u)v_{\Cal K}(D(\Cal F))$$

This  implies 
$$\pr _1(j(M_u/K_{\bar u^n}))\leqslant 2e_{1u}v^1_{\Cal K}(D(\Cal F))<
e_{1u}(p^{u_N^0}c^*_1-1)\leqslant e_{1u}p^uc^*_1-1\leqslant\pr _1(j_u)$$
The lemma is proved.
\enddemo 

Continue the proof of our proposition.

The above lemma implies the existence of a unique root $\theta _u$ 
of $F_u$ such that 
$$v_{M_u}(\theta _{u+1}^p-\theta _u)=i_u+(0,\dots ,0,1)$$
where $\varphi _{M_{\bar u^n}/K_{\bar u^n}}(i_u)=j_u$. 
Similarly to the proof of Lemma 6.2    
$$\frac{j_u-j(M_u/K_{\bar u^n})}{i_u-i(M_u/K_{\bar u^n})}=
\bar e^{-1}(M_u/K_{\bar u^n})$$
where $\bar e(M_u/K_{\bar u^n})=(1,\dots ,1,d)$. Applying Corollary 3.4    
we obtain 
$$i_u=(1,\dots ,1,d)j_u-v_{K_{\bar u^n}}(D(F_u))+(0,\dots ,0,d-1).$$

Therefore, 
$$\pr _1(i_u)+\delta _{1N}\geqslant\pr _1j_u+\delta _{1N}-v^1_{K_{\bar u^n}}(D(F_u))
\geqslant e_{1u}p^uc^*_1-e_{1u}v^1_{\Cal K}(D(F))
\geqslant e_{1u}p^uc^*_1/2$$
i.e. $\eta _u-\eta _{u+1}^p\in \m ^1_{M_u}(p^uc^*_1/2)$. 

As earlier, the uniqueness of $\theta _u$ implies the existence of 
$\theta _{u+1}:=\eta _{u+1}$ such that $\theta _u=\eta _u$ and the required 
sequence $\{\eta _u\ |\ u\geqslant u^0_N\}$ is constructed. 

The uniqueness property of $\theta _u=\eta _u$ implies also that 
$K_{\bar u^{n+1}}(\eta _{u+1})=K_{\bar u^{n+1}}(\eta _u)$.

Consider the tower $M_{\centerdot }=M_{u^0_N}K_{\centerdot }$. 
Then 
\medskip 

--- $M_{\centerdot }\in\Cal B^{fa}(N)$ and has parameters 
$\bar u^*(K_{\centerdot })$ and $c^*(u^0_N, K_{\centerdot })/2$; 
\medskip

--- $\Cal K_{\bar u^0}(\eta )\in X(M_{\centerdot }; m_1,\dots ,m_{N-1})$, where  
$\eta =\varprojlim\eta _u$ is a root of $\Cal F(T)$ in $R_0(N)$. 
The choice of this root $\eta $ of $\Cal F(T)$ determines an embedding of $\Cal L$ 
into $R_0(N)$ which induces the identity on $\Cal K$;
\medskip 

--- by taking $L_{\centerdot }=K_{u^0_1}(\eta _{u^0_N})K_{\centerdot }$, we obtain 
a separable extension of $K_{\centerdot }$ with the parameters $u^0_1$ and $c^*(u^0_N,K_{\centerdot })/2$ 
such that 
$\Cal L\in X(L_{\centerdot })$.

The proposition is proved. 
\enddemo 

\proclaim{Corollary 6.6} Suppose $K_{\centerdot }\in\Cal B^{fa}(N)$ with parameters 
$\bar u^0(K_{\centerdot })=(u^0_1,\dots ,u^0_N)$ and $c^*(u^0_N,K_{\centerdot })$. 
Suppose that 
$\Cal K\in X(K_{\centerdot })$ and 
$\Cal L/\Cal K$ is a finite separable extension in $R_0(N)$ with 
standard $F$-structure. Then there is an $L_{\centerdot }\in\Cal B^{fa}(N)$ such that 
\newline 
{\rm a)} $L_{\centerdot }$ is a finite separable extension of $K_{\centerdot }$;
\newline 
{\rm b)} $\Cal L\in X(L_{\centerdot })$; 
\newline 
{\rm c)} $L_{\centerdot }$ has parameters of the form $v^0=v^0(L_{\centerdot })$ 
and $c^*(v^0,L_{\centerdot })=c^*(u^0_N, K_{\centerdot })/2^N$. 
\endproclaim 

\demo{Proof} Apply the construction from the proof of the above 
proposition to the sequence of extensions 
$$\Cal K\subset \Cal K\Cal L(1)\subset\dots \subset \Cal K\Cal L(N)=\Cal L$$
and notice that for sufficiently large first parameters $n^*_i=n^*((L_{\centerdot}(i))$, 
where $1\leqslant i\leqslant N$, the second parameter can be taken in the form   
$c^*(n^*_i, L_{\centerdot }(i))=c^*(u^0_N, K_{\centerdot })/2^i$. 
\enddemo 

\proclaim{Corollary 6.7} The correspondence 
$K_{\centerdot }\mapsto\cl (K)\in\widetilde{\LF }_R(N)$, where 
$K_{\centerdot }\in\Cal B^{fa}(N)$, induces the identification 
of absolute Galois groups $\psi :\Gamma _{\widetilde{K}}\longrightarrow\Gamma _{\Cal K}$
(here $\widetilde{K}$ is the $p$-adic closure of the $\cup_{m\geqslant 0}K_m$). This identification 
is compatible with ramification filtrations, i.e. 
for any $j\in J(N)$,  
$$\Gamma _{K_0}^{(j)}\cap\Gamma _{\widetilde{K}}=\Gamma _{\Cal K}^{(\varphi _{\widetilde{K}/K_0}(j))}$$
(where $\varphi _{\widetilde{K}/K_0}=\lim _{m\to\infty }\varphi _{K_m/K_0}$). 
\endproclaim

6.3. The above results give that if $\Cal K\in X(K_{\centerdot })$ 
with $K_{\centerdot }\in\Cal B^{fa}(N)$, then $R_0(N)$ contains 
a separable closure of $\Cal K$. Because $R_0(N)$ is perfect 
the algebraic closure of $\Cal K$ in $R_0(N)$ is algebraically closed. 
Even more, 
$R_0(N)$ is $v^1_{\Cal K}$-complete, therefore,  
$R_0(N)$ contains the $v^1_{\Cal K}$-completion $\Cal R(\bar\Cal K)$ 
of $\bar\Cal K$. 

\proclaim{Proposition 6.8} $R_0(N)=\Cal R(\bar\Cal K)$.
\endproclaim 

\demo{Proof} Suppose $K_{\centerdot }$ has parameters $n^*=n^*(K_{\centerdot })$ 
and $c^*=c^*(n^*, K_{\centerdot })$. 
Let $c^{**}=c^*(n^*, K_{\centerdot })/2^N$. The proposition easily follows from 
the following Lemma.

\proclaim{Lemma 6.9}  For any $\alpha\in R(N)$, there is a 
finite separable extension $\Cal L$ of $\Cal K$ and $\beta\in\Cal O_{\Cal R(\Cal L)}$ 
such that $\alpha\equiv\beta\operatorname{mod}\m ^1_{\Cal K}(c^{**})$. 
\endproclaim 

\demo{Proof of Lemma} Suppose $\alpha =(a_u\operatorname{mod}p)_{u\geqslant 0}$, where 
$a_u\in\Cal O_{\Bbb C(N)_p}$ and $a_{u+1}^p\equiv a_u\operatorname{mod}p$ for all $u\geqslant 0$. 
We can assume that $a_0\in L_0$, where $L_0$ is a finite 
extension of $K_0$ such that $L_{\centerdot }=L_0K_{\centerdot }\in\Cal B^{fa}(N)$. 

By the above Corollary 6.6,  $L_{\centerdot }$ has parameters 
$m^*=m^*(L_{\centerdot })\geqslant n^*$ and $c^{**}$. 

Suppose $\bar v^0=(v^0_1,\dots ,v^0_N)$ with $v^0_N\geqslant m^*$  
 is an index parameter from the 
construction of some $\Cal L'\in X(L_{\centerdot })$.

Then $\Cal O_{L_0}\subset \Cal O_{L_{\bar v^0}}$ 
and  
$\pr _{v^0_N}(\Cal O_{\Cal L'})=\Cal O_{L_{\bar v^0}}\operatorname{mod}\m ^1_{\Cal K}(c^{**})$. 

In particular, there is an $\alpha '\in\Cal O_{\Cal L'}$ such that 
$\pr _{v^0_N}\alpha '=a_0\operatorname{mod}\m ^1_{K_{\centerdot }}(c^{**})$, 
or equivalently, 
$$\sigma ^{v^0_N}\alpha \equiv\alpha '\operatorname{mod}\m ^1_{\Cal K}(p^{v^0_N}c^{**}).$$
Therefore, $\alpha\equiv \sigma ^{-v^0_N}\alpha '\operatorname{mod}\m ^1_{\Cal K}(c^{**})$, 
and the lemma is proved because $\sigma ^{-v^0_N}\alpha '\in\Cal O_{\Cal R(\Cal L)}$, 
where $\Cal L$ is a separable extension of $\Cal K$ such that $\cl (\Cal L')=\cl (\Cal L)$. 
\enddemo 
\enddemo

\subhead 7. The functors $\Cal X$ and $\Cal X_{K_{\centerdot }}$ 
\endsubhead

\subsubhead {\rm 7.1.} The functor $\Cal X_{K_{\centerdot }}$, $K_{\centerdot }\in\Cal B^{fa}(N)$
\endsubsubhead 

Let $K_{\centerdot }\in\Cal B^{fa}(N)$ and let $\Cal B^{a}_{K_{\centerdot }}(N)$ be the 
category of finite separable extensions $L_{\centerdot }$ of $K_{\centerdot }$ in $\Cal B^a(N)$. 
Morphisms in $\Cal B_{K_{\centerdot }}^{a}(N)$ are morphisms $f$ in the category $\Cal B^a(N)$ such that 
$f|_{K_{\centerdot }}=\id $. 

Let $\widetilde{\LF }_R(N)_{K_{\centerdot }}$ be the category of finite 
separable extensions of $\cl (\Cal K)\in\widetilde{\LF }_R(N)$, where 
$\Cal K\in X(K_{\centerdot })$. In this section we use  results of n.6 
about the correspondence $K_{\centerdot }\mapsto \cl (\Cal K)$, where 
$K_{\centerdot }\in\Cal B^{fa}(N)$, to construct an equivalence between the 
categories $\Cal B_{K_{\centerdot }}^{a}(N)$ and $\widetilde{\LF }_R(N)_{K_{\centerdot }}$.
\medskip 

Let $L_{\centerdot }$ be a finite separable extension of $K_{\centerdot }\in\Cal B^{fa}(N)$ 
in $\Cal B(N)$. 
Then $L_{\centerdot }\in\Cal B^a(N)$, cf. n.4.2.  
Choose a finite Galois extension $E_{\centerdot }$ of $K_{\centerdot }$ such that 
$E_{\centerdot }\in\Cal B^{fa}(N)$ and $E_{\centerdot }\supset L_{\centerdot }$, cf. Prop.4.4.  
If $\Cal K\in X(K_{\centerdot })$ then there is a unique separable extension $\Cal E$ 
of $\Cal K$ in $R_0(N)$ such that $\Cal E\in X(E_{\centerdot })$ and 
$[\Cal E:\Cal K]=[E_{\centerdot }:K_{\centerdot }]$, cf. n.6. 
Therefore, $G=\Gal (E_{\centerdot }/L_{\centerdot })$ acts on $\Cal E$ and we can set 
$\Cal L=\Cal E^H$, where $H\subset G$ is such that $E_{\centerdot }^H=L_{\centerdot }$. 

\proclaim{Proposition 7.1} With the above notation, $\cl (\Cal L)\in\widetilde{\LF}_R(N)$ does not 
depend on the choice of $\Cal K \in X(K_{\centerdot })$ and 
$E_{\centerdot }\in\Cal B^{fa}(N)$.
\endproclaim 

The proof is straightforward. 
\medskip 

Suppose $L_{\centerdot }, L_{\centerdot }'\in\Cal B^a_{K_{\centerdot }}(N)$ 
and $f:L_{\centerdot }\longrightarrow L_{\centerdot }'$ is a morphism 
in $\Cal B^a_{K_{\centerdot }}(N)$. In other words, 
$f$ is a $P$-continuous and compatible with the corresponding 
$F$-structures automorphism of $\Bbb C(N)_p$ 
such that $f(L_m)=L_m'$ for $m\gg 0$ and $f|_{K_{\centerdot }}=\id $. 

Choose $E_{\centerdot }\in\Cal B^{fa}(N)$ such that $E_{\centerdot }\supset L_{\centerdot }$ and 
$E_{\centerdot }$ is finite Galois over $K_{\centerdot }$. Let $E'_{\centerdot }=f(E_{\centerdot })$, 
$G=\Gal (E_{\centerdot }/K_{\centerdot })$, 
$G'=\Gal (E'_{\centerdot }/K_{\centerdot })$, $H=\Gal (E_{\centerdot }/L_{\centerdot })$ and 
$H'=\Gal (E'_{\centerdot }/L_{\centerdot }')$. 

Let $\Cal K\in X(K_{\centerdot })$ and let $\Cal E$ be its field Galois extension in $X(E_{\centerdot })$ 
of degree $[E_{\centerdot }:K_{\centerdot }]$. Let $f_R$ be an automorphism of $R_0(N)$ 
induced by $f$. Then $f_R$ is $P$-continuous and compatible with $F$-structures, 
$f_R(\Cal E)\in X(E'_{\centerdot })$ and $f_R(\Cal E)^{H'}=f(\Cal E^H)\in X(L_{\centerdot }')$. 

So, $f_R\in\Hom _{\widetilde{\LF }_R(N)_{K_{\centerdot }}}(\Cal X_{K_{\centerdot }}(L_{\centerdot }),
\Cal X_{K_{\centerdot }}(L_{\centerdot }'))$. Clearly, if we set $f_R=\Cal X_{K_{\centerdot }}(f)$ then we get 
a functor $\Cal X_{K_{\centerdot }}$ from $\Cal B_{K_{\centerdot }}^{a}(N)$ to 
$\widetilde{\LF }_R(N)_{K_{\centerdot }}$. 
\medskip 

Summarizing the results of n.6 we obtain the following principal result of this paper.   

\proclaim{Theorem 2}  {\rm a)} The above defined functor $\Cal X_{K_{\centerdot }}$, 
where $K_{\centerdot }\in\Cal B^{fa}(N)$, 
is an equivalence of the categories $\Cal B^{a}_{K_{\centerdot }}(N)$ and 
$\widetilde{\LF }_R(N)_{K_{\centerdot }}$; 
\newline 
{\rm b)} $\Cal X_{K_{\centerdot }}$ induces an identification 
$\psi _{K_{\centerdot }}$ of groups 
$\Gamma _{\widetilde{K}}=\Gal (\bar K/\widetilde{K})$ and 
\linebreak 
$\Gamma _{\Cal K}=\Gal (\Cal K_{\sep }/\Cal K)$, where 
$\Cal K\in X(K_{\centerdot })$ and $\Cal K_{\sep }$ is 
the separable closure of $\Cal K$ in $R_0(N)$;
\newline 
{\rm c)} the identification $\psi _{K_{\centerdot }}$ is compatible with ramification 
filtrations on $\Gamma _{\widetilde{K}}$ and $\Gamma _{\Cal K}$, i.e. 
for any $j\in J(N)$, $\psi _{K_{\centerdot }}$ identifies  the groups 
$\Gamma _{\widetilde{K}}\cap\Gamma _{K_0}^{\varphi _{\widetilde{K}/K_0}(j)}$ and 
$\Gamma _{\Cal K}^{(j)}$.
\endproclaim 
\medskip 

\subsubhead {\rm 7.2} The functor $\Cal X:\Cal B^a(N)\longrightarrow \RLF _{\Cal R}(N)$ 
\endsubsubhead 

Let $\RLF _{R}(N)$ be the category of $P$-closed perfect subfields in  
$R_0(N)$. These subfields are considered with 
their natural $F$-structure and $P$-topology. Morphisms are $P$-continuous isomorphisms 
of such fields, which are compatible with corresponding $F$-structures. 

If $K_{\centerdot }\in\Cal B^a(N)$ choose $L_{\centerdot }\in\Cal B^{fa}(N)$ such that 
$L_{\centerdot }/K_{\centerdot }$ is a finite Galois extension. 
If $\Cal L\in X(L_{\centerdot })$ then $G=\Gal (L_{\centerdot }/K_{\centerdot })$ acts on 
$\Cal R(\Cal L)$. Indeed, 
for any $g\in G$ and any $\Cal L\in X(L_{\centerdot })$, the action 
of $g$ on $L_{\centerdot }$ induces a field isomorphism 
$g:\Cal L\longrightarrow\Cal L'$, where 
$\Cal L'\in X(L_{\centerdot })$ and we have a natural identification 
$\Cal R(\Cal L)=\Cal R(\Cal L')$. With the above notation 
set $\Cal X(K_{\centerdot })=\Cal R(\Cal L)^G\in\RLF _R(N)$. 

\proclaim{Proposition 7.2} $\Cal X(K_{\centerdot })$ does not depend on a choice 
of $L_{\centerdot }\in\Cal B^{fa}(N)$. 
\endproclaim 

\demo{Proof} Suppose $L_{\centerdot }'\in\Cal B^{fa}(N)$ is such that 
$L_{\centerdot }'/K_{\centerdot }$ is a finite Galois extension with the Galois group 
$G'$. Choose $M_{\centerdot }\in\Cal B^{fa}(N)$ such that $M_{\centerdot }\supset L_{\centerdot }$, 
$M_{\centerdot }\supset L_{\centerdot }'$ and $M_{\centerdot }$ is a finite Galois 
extension of $K_{\centerdot }$ with the Galois group $S$. 

Let $H=\Gal (M_{\centerdot }/L_{\centerdot })$, 
$H'=\Gal (M_{\centerdot }/L_{\centerdot }')$. If $\Cal L\in X(L_{\centerdot })$, 
$\Cal L'\in X(L_{\centerdot }')$, then there are $\Cal M\in X(M_{\centerdot })$ and 
$\Cal M'\in X(M_{\centerdot })$ such that $\Cal M/\Cal L$ and $\Cal M'/\Cal L'$ 
are Galois extensions with Galois groups $H$ and $H'$, respectively. 
Then $\Cal R(\Cal M)=\Cal R(\Cal M')$ and 
$\Cal R(\Cal L')^{G'}=\Cal R(\Cal M')^S=\Cal R(\Cal M)^S=\Cal R(\Cal L)^G$.

The proposition is proved.
\enddemo 

Suppose $K_{\centerdot }, K_{\centerdot }'\in\Cal B^a(N)$ and 
$f\in\Hom _{\Cal B^a(N)}(K_{\centerdot }, K_{\centerdot }')$, i.e. 
$f:\Bbb C(N)_p\longrightarrow\Bbb C(N)_p$ is a $P$-continuous 
and compatible with $F$-structure 
field automorphism such that $f(K_{\centerdot })=K_{\centerdot }'$. As earlier, 
denote by $f_R$ the automorphism of $R_0(N)$ which is induced by $f$. 

Choose $L_{\centerdot }\in\Cal B^{fa}(N)$ such that $L_{\centerdot }/K_{\centerdot }$ 
is a finite Galois extension with the group $G$. Then $L_{\centerdot }'=f(L_{\centerdot })$ 
is a Galois extension of $K_{\centerdot }'$ with the same group $G$. 
If $\Cal L\in X(L_{\centerdot })$ then $f_R(\Cal L)=\Cal L'\in X(L_{\centerdot }')$ and 
$f_R(\Cal X(K_{\centerdot }))=f_R(\Cal R(\Cal L)^G)=\Cal R(\Cal L')^G=\Cal X(K_{\centerdot }')$. 

So, $f_R\in\Hom _{\RLF _R(N)}(\Cal X(K_{\centerdot }),\Cal X(K_{\centerdot }'))$ 
and $\Cal X:\Cal B^a(N)\longrightarrow\RLF _R(N)$ is a functor. The following property 
follows directly from the above definitions. 

\proclaim{Proposition 7.3} {\rm a)} $\Cal X$ is a strict functor;
\newline
{\rm b)} if $L_{\centerdot },K_{\centerdot }\in\Cal B^{fa}(N)$ and 
$L_{\centerdot }$ is a finite separable extension of 
$K_{\centerdot }$ then 
$$\Cal R(\Cal X_{K_{\centerdot }}(L_{\centerdot }))=\Cal X(L_{\centerdot }).$$
\endproclaim 

7.3. Let $\varepsilon =(\varepsilon ^{(n)}\operatorname{mod}p)_{n\geqslant 0}\in R(1)\subset R(N)$, 
where $\varepsilon ^{(0)}=1$, $\varepsilon ^{(1)}\ne 1$ and ${\varepsilon ^{(n+1)}}^p=\varepsilon ^{(n)}$ 
for all $n\geqslant 0$, be Fontaine's element. 
Let $<\varepsilon>=\varepsilon ^{\Bbb Z_p}\subset R(1)^*$ be the multiplicative 
subgroup of all Fontaine's elements. Notice, if 
$f:\Bbb C(N)_p\longrightarrow\Bbb C(N)_p$ is 
a field automorphism then  
$f_R(<\varepsilon >)=<\varepsilon >$, where $f_R$ is induced by $f$.

\proclaim{Lemma 7.4} The correspondence $f\mapsto f_R$ identifies  
$\Aut \Bbb C(N)_p$ and the subgroup 
$\Aut 'R_0(N)$ of $g\in\Aut R_0(N)$ such that $g(<\varepsilon >)=<\varepsilon >$.
\endproclaim  

\demo{Proof} We have noticed already that for any $f\in\Aut \Bbb C(N)_p$, $f_R(<\varepsilon >)=<\varepsilon >$. 

Suppose $g\in\Aut R_0(N)$ and $g(<\varepsilon >)=<\varepsilon >$, 
i.e. $g(\varepsilon )=\varepsilon ^a$ with $a\in\Bbb Z_p^*$. 

Notice that $g:R(N)\longrightarrow R(N)$ induces 
the automorphism $W(g):W(R(N))\longrightarrow W(R(N))$, where $W$ is the functor of Witt vectors.  
Consider the Fontaine map 
$$\gamma :W(R(N))\longrightarrow \Cal O_{\Bbb C(N)_p}$$ 
given by the correspondence $(r_0,r_1,\dots ,r_n,\dots )\mapsto 
r^{(0)}+pr^{(1)}+\dots +p^nr^{(n)}+\dots $, 
where for any $r=(r_m\operatorname{mod}p)_{m\geqslant 0}\in R(N)$ and $n\geqslant 0$, 
$r^{(n)}=\lim _{m\to\infty }r_{m+n}^m$. This map is a surjective 
morphism of $p$-adic algebras and its kernel $J$ is a principal ideal 
generated by $1+[\varepsilon ]^{1/p}+\dots +[\varepsilon ]^{(p-1)/p}$. 
Therefore, $W(g)(J)=J$ and $W(g)$ induces an  
automorphism $f=W(g)\operatorname{mod}J$ of $\Bbb C(N)_p$. Clearly, $f_R=g$. 
The lemma is proved. 
\enddemo 

\remark{Remark} From the above description of the correspondence 
$f\mapsto f_R$ it easily follows that 
$f$ is $P$-continuous (resp., compatible with $F$-structure) if and only if 
$f_R$ possess the same property. 
\endremark 
\medskip 

7.4. Introduce the following definition.  

\definition{Definition} A subfield $\widetilde{K}$ of $\Bbb C(N)_p$ is an SAPF-field if 
there is $K_{\centerdot }\in\Cal B^a(N)$ such that $\widetilde{K}$ is the $p$-adic closure 
of $\cup _{n\geqslant 0}K_n$. 
\enddefinition 

\remark{Remark} The above definied $\SAPF$-fields are higher dimensional analogues of 
strict arithmetic profinite extensions introduced in [FW1-2]. 
\endremark 

Denote by $\SAPF (N)$ the category of $\SAPF$-fields in $\Bbb C(N)_p$, such that 
if $\widetilde{K},\widetilde{K}'\in\SAPF (N)$, then 
$\Hom _{\SAPF(N)}(\widetilde{K},\widetilde{K}')$ consists of 
$P$-continuous and compatible with $F$-structures 
$f\in\Aut \Bbb C(N)_p$ such that $f(\widetilde{K})=\widetilde{K}'$. 

Let $\widetilde{K}\in\SAPF (N)$. Set $\widetilde{\Cal X}(\widetilde{K})=
\Cal X(K_{\centerdot })$, where $K_{\centerdot }\in\Cal B^a(N)$ is such that 
$\widetilde{K}$ is a $p$-adic closure of $\cup _{n\geqslant 0}K_n$. 

\proclaim{Lemma 7.5} The above defined $\widetilde{\Cal X}(\widetilde{K})$ 
does not depend on the choice of $K_{\centerdot }\in\Cal B^a(N)$.
\endproclaim 

\demo{Proof} The proof follows directly from the construction of the functor 
$\Cal X$ and proposition 5.6. 
\enddemo 

The correspondence $\widetilde{K}\mapsto\widetilde{\Cal X}(\widetilde{K})$ is naturally 
extended to the functor 
\linebreak 
$\widetilde{\Cal X}:\SAPF (N)\longrightarrow\RLF _R(N)$.

Taking together the above results about the functor $\Cal X$ we obtain the following theorem.  

\proclaim{Theorem 3} Suppose $K_{\centerdot }\in\Cal B^a(N)$ and 
$\widetilde{K}$ is a $p$-adic closure of $\cup _{n\geqslant 0}K_n$.  
Then the functor $\widetilde{\Cal X}$ induces the identification 
$\iota :\Gamma _{\widetilde{K}}\longrightarrow \Gamma _{\widetilde{\Cal K}}$ 
where $\widetilde{\Cal K}=\Cal X(K_{\centerdot })$. If $K_{\centerdot }\in\Cal B^{fa}(N)$ 
and $\Cal K\in X(K_{\centerdot })$ then $\Cal R(\Cal K)=\widetilde{\Cal K}$ and 
under a natural identification $\Gamma _{\Cal K}=\Gamma _{\widetilde{\Cal K}}$, 
the identification 
$\iota $ is compatible with ramification filtrations, i.e. for any 
$j\in J(N)$, 
$$\Gamma _{\widetilde{K}}\cap\Gamma _{K_0}^{(\varphi _{\widetilde{K}/K_0}(j)}
=\Gamma _{\Cal K}^{(j)}.$$
\endproclaim 
\medskip 

\subhead 8. A property of the $P$-continuity for the functor $\Cal X$ 
\endsubhead 
\medskip 

8.1. Suppose $\Cal K\in\LF _p(N)$. 

Let $\Gamma _{\Cal K}^{\ab }(p)$ be the Galois group 
of the maximal abelian $p$-extension of $\Cal K$. 

For any $M\geqslant 1$, consider the Witt-Artin-Schreier duality
$$\Gamma _{\Cal K}^{\ab }(p)/p^M \times W_M(\Cal K)/(\sigma -\id )W_M(\Cal K)
\longrightarrow W_M(\Bbb F_p)$$
where $\sigma $ is the Frobenius endomorphism of the additive group $W_M(\Cal K)$ 
of Witt vectors of length $M$ with coefficients in $\Cal K$. This allows us to provide  
$\Gamma _{\Cal K}^{\ab }(p)/p^M$ with 
the $P$-topological structure. 
Its basis of open 0-neighborhoods consists of the annihilators of the  
compact subsets of $W_M(\Cal K)/(\sigma -\id )W_M(\Cal K)$. 
By results of n.1.2 the basis of such compact subsets consists 
of the images in $W_M(\Cal K)/(\sigma -\id )W_M(\Cal K)$ of all 
subsets of the form 
$$W_M(D)=\{(a_0,\dots ,a_{M-1})
\in W_M(\Cal K)\ |\ a_0,\dots ,a_{M-1}\in D\}$$ 
where $D\in\Cal C_b(\Cal K)$ is the basis of compact subsets in $\Cal K$.   

Finally, the $P$-topology on 
$\Gamma _{\Cal K}^{\ab }(p)$ appears as the projective limit topology 
of the projective system of $P$-topological groups $\Gamma ^{\ab }_{\Cal K}(p)/p^M$. 
\medskip 

8.2. Suppose $K\in\LF _0(N)$ and $K$ contains a primitive 
$p^M$-th root of unity $\zeta _{p^M}$. Then the $P$-topological 
structure on $K^*$ induces the $P$-topological structure on 
$\Gamma _K^{\ab }(p)/p^M$, where $\Gamma _K(p)$ is the 
Galois group of the maximal abelian $p$-extension of $K$. 
This structure is defined similarly to the characteristic $p$ case 
by the use of the Kummer duality 

$$\Gamma _K^{\ab }(p)/p^M \times K^*/{K^*}^{p^M}\longrightarrow <\zeta _{p^M}>$$

We don't need this structure in a full generality. 
Let $\widetilde{\Gamma }_K(p)/p^M$ be the quotient 
of $\Gamma _K^{\ab }(p)/p^M$ by the annihilator of the subgroup 
in $K^*$ generated by the elements of $(1+p\Cal O_K)^*$. 
Then we have the induced pairing 
$$\widetilde{\Gamma }_K^{\ab }(p)/p^M\times (1+p\Cal O_K)^*\longrightarrow <\zeta _{p^M}>$$
and a basis of open subgroups in $\widetilde{\Gamma }_K^{\ab }(p)/p^M$ consists of the 
annihilators of the subsets  $1+pD$, where $D\in\Cal C_b(K)$ and 
$\Cal C_b(K)$ is a basis of compact subsets in $K$ from n.1.2.  
\medskip 

8.3. Suppose $K_{\centerdot }\in\Cal B^{fa}(N)$ and for a sufficiently 
large $n$, $K_n$ contains a primitive $p^M$-th root of unity. 

Let $\widetilde{K}$ be the $p$-adic closure of $\cup _{n\geqslant 0}K_n$. 
Then for any $M\in\Bbb N$, we have a natural identification 
$$\Gamma ^{\ab }_{\widetilde{K}}(p)/p^M=\mathbin{\underset{n}\to\varprojlim}
\Gamma _{K_n}^{\ab }/p^M.$$
Applying the arguments from n.2.9 we can also write 
$$\Gamma ^{\ab }_{\widetilde{K}}(p)/p^M=\mathbin{\underset{n}\to\varprojlim}
\widetilde{\Gamma }_{K_n}^{\ab }/p^M.$$
Therefore, the basis of $P$-open neighborhoods in $\Gamma _{\widetilde{K}}^{\ab }(p)/p^M$ 
consists of annihilators of all compact subsets 
$1+pD\subset (1+p\Cal O_{\widetilde{K}})^{\times }$, 
where $D\in\Cal C_b(K_n)$ for some $n\geqslant 0$. 
\medskip 

8.4. Suppose $K_{\centerdot }\in\Cal B^{fa}(N)$, 
$\Cal K\in X(K_{\centerdot })$ and 
$\iota :\Gamma _{\widetilde{K}}\longrightarrow\Gamma _{\Cal K}$ 
is the identification of 
Galois groups (where $\widetilde{K}$ is the $p$-adic closure 
of the $\cup _{n\geqslant 0}K_n$) from Theorem 3. 
Suppose for each $M\in\Bbb N$, $\zeta _{p^M}\in K_n$ if $n\gg 0$ 
and consider the groups 
$\Gamma _{\widetilde{K}}^{\ab }/p^M=
\mathbin{\underset{n}\to\varprojlim}\widetilde{\Gamma }_{K_n}/p^M$ 
and $\Gamma _{\Cal K}^{\ab }/p^M$ with the above $P$-topological 
structures. 

\proclaim{Theorem 4} With the above notation, the identification 
$$\iota\operatorname{mod}p^M:\Gamma ^{\ab }_{\widetilde{K}}/p^M
\longrightarrow\Gamma _{\Cal K}^{\ab }/p^M$$
is $P$-continuous.
\endproclaim 

\demo{Proof} 

8.4.1. Consider the dual morphism 
$$\tilde \iota _M:W_M(\Cal K)/(\sigma -\id )W_M(\Cal K)\longrightarrow 
\widetilde{K}^*/\widetilde{K}^{*p^M}.$$ 
Then $\iota \operatorname{mod}p^M$ is $P$-continuous if and only if 
$\tilde\iota _M$ transforms each P-compact subset in 
$W_M(\Cal K)/(\sigma -\id )W_M(\Cal K)$ onto a $P$-compact subset in 
$\widetilde{K}^*/\widetilde{K}^{p^M}$. 

Notice that the map $\tilde\iota _M$ can be characterised as follows. 

Choose a primitive $p^M$-th root of unity $\zeta _M$. 
Let $\bar w\in W_M(\Cal K)/(\sigma -\id )W_M(\Cal K)$ and let 
$w\in W_M(\Cal K)$ be a lifting of $\bar w$. Consider $T\in W_M(R(N))$ such that 
$\sigma T-T=w$ 
then for any $\tau\in\Gamma _{\Cal K}$, $\tau T-T=a_{\tau }\in W_M(\Bbb F_p)$. 
Let $\bar v\in \widetilde{K}^*/\widetilde{K}^{*p^M}$ and $v\in \widetilde{K}^*$ be a lifting of $\bar v$. 
Consider $Z\in\Bbb C(N)_p$ such that $Z^{p^M}=v$. Then 
for any $\tau\in\Gamma _{\widetilde{K}}$, $\tau Z/Z=\zeta _M^{b_\tau }$, 
where $b_{\tau}\in W_M(\Bbb F_p)$. 
With the above notation, with respect to the identification 
$\Gamma _{\widetilde{K}}=\Gamma _{\Cal K}$ given by the construction 
of the functor $\Cal X_{K_{\centerdot }}$, we have the following criterion:

$$\tilde\iota _M(\bar w)=\bar v\ \ \ \Leftrightarrow\ \ \ 
a_{\tau }=b_{\tau }\ \ \forall\tau\in \Gamma _{\widetilde{K}}=\Gamma _{\Cal K}$$
\medskip

8.4.2. 
As earlier, let $\Cal R(\Cal K)$ be the completion of the radical closure of $\Cal K$ 
(with respect to 1st valuation). Denote by 
$\Cal R(\Cal O_{\Cal K})$ its valuation ring.  

Notice first that the natural embedding $\Cal K\subset\Cal R(\Cal K)$ induces 
a natural identification of $P$-topological groups 
$$W_M(\Cal R(\Cal K))/(\sigma -\id )W_M(\Cal R(\Cal K))=
W_M(\Cal K)/(\sigma -\id )W_M(\Cal K)$$

Let $\varepsilon $ be Fontaine's elements. Recall, 
$\varepsilon =(\varepsilon ^{(n)})_{n\geqslant 0}\in R=R(1)\subset R(N)$ is such that 
$\varepsilon ^{(0)}=1$, $\varepsilon ^{(1)}\ne 1$ and  we can assume that 
$\varepsilon ^{(M)}=\zeta _M$ --- this is the primitive $p^M$-th root of unity chosen in 8.4.1.  
From the construction of $\Cal K\in X(K_{\centerdot })$ it follows that 
$\varepsilon\in\Cal R(\Cal O_{\Cal K})$. Consider the map 
$$\pr :\frac{1}{[\varepsilon ]-1}W_M(\Cal R(\Cal O_{\Cal K}))
\longrightarrow W_M(\Cal R(\Cal K))/(\sigma -\id )W_M(\Cal R(\Cal K))$$
induced by the projection $W_M(\Cal R(\Cal K))
\longrightarrow W_M(\Cal R(\Cal K))/(\sigma -\id )W_M(\Cal R(\Cal K))$. 

\proclaim{Lemma 8.1} $\pr $ is surjective.
\endproclaim  

\demo{Proof} This follows from the formula 
$$W_M(\Cal K)=\bigcup\Sb s\geqslant 0\endSb \frac{1}{\sigma ^s([\varepsilon ]-1)}W_M(\Cal O_{\Cal K}).$$
\enddemo 

\remark{Remark} It can be easily seen that the family of sets 
$$\left\{ \pr \left (\frac{1}{[\varepsilon ]-1}
W_M(\sigma ^{-s}D)\right )\ |\ s\in\Bbb Z_{\geqslant 0}, D\subset \Cal O_{\Cal K}, D\in\Cal C_b(\Cal K)\right\}$$
is a basis of compact subsets in $W_M(\Cal K)/(\sigma -\id )W_M(\Cal K)$.
\endremark 

Let $w\in W_M(\Cal K)$ be the element from n.8.4.1. By the above lemma, 
there is an $f\in W(\Cal R(\Cal O_{\Cal K}))$ such that 
$w=f/([\varepsilon ]-1)\operatorname{mod}p^M$. Therefore, 
if $U\in W(R_0(N))$ is such that 
$\sigma U-U=f/([\varepsilon ]-1)$ then for any $\tau\in\Gamma _{\Cal K}$, 
$\tau U-U=\tilde a_{\tau }\in W(\Bbb F_p)$, where 
$\tilde a_{\tau }\operatorname{mod}p^M=a_{\tau }$. 
\medskip 

8.4.3. Let $\varepsilon _1=\sigma ^{-1}\varepsilon $, then 
$$s=([\varepsilon ]-1)/([\varepsilon _1]-1)\in W^1(R(1))\subset W(R(1))\subset W(R(N))$$ 
where $W^1(R(1))=\Ker \gamma :W(R(1))\longrightarrow\Cal O_{\Bbb C_p}$ 
is Fontaine's map. It is known [Ab2] that $s$ generates the ideal $W^1(R(1))$. 
Notice that similar arguments show that 
$s$ generates also the kernel $W^1(R(N))$ of the analogue of Fontaine's map from 
$W(R(N))$ to $\Cal O_{\Bbb C(N)_p}$.

Let $T_1=U([\varepsilon _1]-1)$. Then $T_1\in W(R(N))$ and 
$\sigma T_1-sT_1=f$. Let $X=U([\varepsilon ]-1)=sT_1\in W^1(R)$, then 
$$\frac{\sigma X}{\sigma s}-X=f$$
and for any $\tau\in\Gamma _{\Cal K}$, $\tau X-X=\tilde a_{\tau }([\varepsilon ]-1)$. 
\medskip 

8.4.4. Let $A(N)_{cris }$ be an analogue of Fontaine's $A_{cris }$ 
constructed by the use of $R(N)$ instead of $R$. This is 
the divided power envelope of the $W(R(N))$ with respect to the 
ideal $W^1(R(N))$, which is generated by $s$. Proceeding as in [Ab2] 
we obtain that if 
$$\frac{\sigma m}{p}-m=f\tag {5}$$
where $m\in \Fil ^1A(N)_{cris }$, then for any 
$\tau\in\Gamma _{\widetilde{K}}$, $\tau m-m=\tilde a_{\tau }\log [\varepsilon ]$. 

Multiplying both parts of the equality (5) by $p$ and taking exponentials  
we obtain the equality 
$$\sigma Y=Y^p\exp (pf)\tag{6}$$
where $Y\in 1+\Fil ^1A(N)_{cris }$ and for any 
$\tau\in\Gamma _{\widetilde{K}}$, $\tau Y/Y=[\varepsilon ]^{\tilde a_{\tau }}$. 
Proceeding again as in [Ab2] we can prove that $Y\in 1+W^1(R(N))$ 
(and therefore can forget about the cristalline ring $A(N)_{cris }$). 
\medskip 

8.4.5. The equation (6) implies that 
$$\sigma ^MY=Y^{p^M}\exp (p\sigma ^{M-1}f+\dots +p^Mf)$$
and, because $\sigma $ is injective on $W(R(N))$, this gives 
$$Y=(\sigma ^{-M}Y)^{p^M}\exp (p\sigma ^{-1}f+\dots +p^M\sigma ^{-M}f)\tag{7}$$
Notice  that for any $\tau\in\Gamma _{\widetilde{K}}$, 
$\tau (\sigma ^{-M}Y)=(\sigma ^{-M}Y)[\sigma ^{-M}\varepsilon ]^{\tilde a_{\tau }}$. 

Apply Fontaine's map $\gamma :W(R(N))\longrightarrow \Cal O_{\Bbb C(N)_p}$ 
to the both parts of ($7$). 
Notice that $\gamma (Y)=1$, $\gamma (\sigma ^{-M}Y)=Z\in 1+p\Cal O_{\Bbb C(N)_p}$, 
$\gamma ([\sigma ^{-M}\varepsilon ])=\zeta _M$ and 
$\gamma (\sigma ^{-s}f)\in \Cal O_{\widetilde{K}}$ for any $s\in\Bbb Z$.  
This gives 
$$Z^{p^M}=\exp (-p\gamma (\sigma ^{-1}f)\dots -p^M\gamma (\sigma ^{-M}f))
\in 1+p\Cal O_{\widetilde{K}}$$
and for any $\tau \in\Gamma _{\widetilde{K}}$, $\tau Z/Z=\zeta _M^{a_{\tau }}$. 
\medskip 

8.4.6. The above computations imply (with the notation from n.8.4.1) that 
if $\bar w=f/([\varepsilon ]-1)\operatorname{mod}(\sigma -\id )W_M(\Cal R(\Cal K))$ 
with $f\in W(\Cal R(\Cal O_{\Cal K}))$ then 
$\tilde \iota _M(\bar w)=\bar v$, where 
$$\bar v=\exp (-p\gamma (\sigma ^{-1}f)-\dots -
p^M\gamma (\sigma ^{-M}f))\operatorname{mod}\widetilde{K}^{*p^M}.$$
It remains to notice that 
by Proposition 1.2 the correspondence 
$$f\mapsto  exp (-p\gamma (\sigma ^{-1}f)-\dots -
p^M\gamma (\sigma ^{-M}f))\operatorname{mod}\widetilde{K}^{*p^M}$$
maps all $P$-compact subsets in $W_M(\Cal R(\Cal O_{\Cal K}))$ 
to $P$-compact subsets in $1+p\Cal O_{\widetilde{K}}$. 

The theorem is proved. 
\enddemo 

\remark{Remark} The above computations in nn.8.4.3-8.4.6 can be used to deduce 
(in the similar way as in [Ab2]) the explicit formula for Hilbert symbol 
for higher dimensional fields from [Vo]. 
\endremark 
\medskip 

\subhead 9. The Grothendieck Conjecture for higher dimensional local fields 
\endsubhead 
\medskip 

9.1. Suppose $K,K'$ are 1-dimensional local fields 
from the category $\LF (1)=\LF _0(1)\coprod\LF _p(1)$. Then any isomorphism 
$f\in\Hom _{\LF (1)}(K,K')$ is given by an automorphism of $\Bbb C(1)_p$ or $\Cal C(1)_p$  
such that $f(K)=K'$. Therefore, $f$ induces the isomorphism of profinite groups 
$$f^*:\Gamma _{K'}\longrightarrow\Gamma _K$$
such that for any $v\geqslant 0$, $f^*(\Gamma _{K'}^{(v)})=\Gamma _K^{(v)}$. 

The inverse statement was proved in [Mo] in the mixed characteristic case 
 and in [Ab4] if the characteristic 
of the residue fields of $K$ and $K'$ is $\geqslant 3$. It is known as 
a local (1-dimensional) analogue of the  Grothendieck Conjecture and can be stated 
in the following form:

{\it If $\iota :\Gamma _K'\longrightarrow\Gamma _{K}$ is an isomorphism 
of profinite groups such that for any $v\geqslant 0$, 
$\iota (\Gamma _{K'}^{(v)})=\Gamma _K^{(v)}$, then there is 
an $f\in\Hom _{\LF (1)}(K,K')$ such that $\iota =f^*$.}
\medskip 

9.2. Suppose $N\geqslant 1$ and $\Cal K,\Cal K'\in\LF _R(N)$. 
Suppose $f\in\Hom _{\LF _R(N)}(\Cal K,\Cal K')$ is isomorphism. 
In other words, 
$f:R_0(N)\longrightarrow R_0(N)$ is  a $P$-continuous 
and compatible with $F$-structures field automorphism such that 
for all $1\leqslant i\leqslant N$, 
$f(\Cal K(i)R(\Cal K(i-1))=\Cal K'(i)R(\Cal K'(i-1))$. Then 
$f^*:\Gamma _{\Cal K'}\longrightarrow\Gamma _{\Cal K}$ 
is an isomorphism of profinite groups such that 
for any $j\in J(N)$, $f^*(\Gamma _{\Cal K'}^{(j)})=\Gamma _{\Cal K}^{(j)}$. 

In addition, suppose $\Cal E$ is a finite extension of $\Cal K$ in $R_0(N)$ 
and $f(\Cal E)=\Cal E'$. Then $\Cal E'$ is a finite extension of 
$\Cal K'$ such that $f^*(\Gamma _{\Cal E'})=\Gamma _{\Cal E}$. Let $M\in\Bbb N$. 
Consider the induced 
isomorphism of the maximal abelian quotients modulo $p^M$-th powers 
$$f^*_M:\Gamma _{\Cal E'}^{\ab }/p^M\longrightarrow\Gamma _{\Cal E}^{\ab }/p^M .$$
It is dual to the isomorphism of additive groups 
$$f_M:W_M(\Cal E)/(\sigma -\id )W_M(\Cal E)\longrightarrow W_M(\Cal E')/(\sigma -\id )W_M(\Cal E').$$
Clearly, $f_M$ is $P$-continuous and, therefore, maps $P$-compact subsets  
to $P$-compact subsets. This implies that $f^*_M$ is $P$-continuous for all $M\in\Bbb N$.

The inverse statement appears as an analogue of the Grothendieck Conjecture 
for higher dimensional local fields of characteristic $p$. 

\proclaim{Theorem 5} With the above notation suppose that $p\geqslant 3$ and 

$$\iota :\Gamma _{\Cal K'}\longrightarrow\Gamma _{\Cal K}$$ 
is an isomorphism of profinite groups such that 
\newline 
{\rm a)} for any $j\in J(N)$, 
$\iota (\Gamma _{\Cal K'}^{(j)})=\Gamma _{\Cal K}^{(j)}$;
\newline 
{\rm b)} if $\Cal E$ and $\Cal E'$ are finite extensions of 
$\Cal K$ and, resp., $\Cal K'$ in $R_0(N)$ such that 
the both $\Cal E$ and $\Cal E'$
have a standard $F$-structure, then for all $M\geqslant 1$, 
the induced isomophism 
$$\iota _M:\Gamma _{\Cal E'}^{\ab }/p^M\longrightarrow
\Gamma _{\Cal E}^{\ab }/p^M$$
is $P$-continuous. 

Then there is an $f\in\Hom _{\LF _R(N)}(\Cal K,\Cal K')$ such that $f^*=\iota $.
\endproclaim 

This statement was proved in [Ab6] in the case $N=2$. The case of general $N$ 
can be done along the same lines.

\remark{Remark} Actually, in the statement of the main theorem in [Ab6] 
there was no requirement that $\Cal E$ and $\Cal E'$ have standard $F$-structure. 
But in the proof we applied this condition only to fields with standard $F$-structure. 
Also, in [Ab6] there was a requirement about the $P$-continuity of the induced group isomorphism 
$\iota ^{\ab }:\Gamma _{\Cal E'}^{\ab }\longrightarrow\Gamma _{\Cal E}^{\ab }$ but 
again in the proof we applied this property only to the induced 
isomorphism of the Galois groups $\Gamma _{\Cal E'}^{\ab }(p)$ and $\Gamma _{\Cal E}^{\ab }(p)$ 
of the maximal $p$-extensions of $\Cal E'$ and $\Cal E$. 
\endremark 
\medskip 

9.3. Suppose $N\geqslant 1$ and $K,K'\in\LF _0(N)$. Any 
$P$-continuous and compatible with $F$-structures field automorphism 
$f:\Bbb C(N)_p\longrightarrow \Bbb C(N)_p$ such that 
$f(K)=K'$ induces an isomphism 
of profinite groups $f^*:\Gamma _{K'}\longrightarrow\Gamma _K$ 
such that $f^*(\Gamma ^{(j)}_{K'})=\Gamma _K^{(j)}$ for any $j\in J(N)$. 

Suppose $E$ is a finite extension of $K$ , 
then $E'=f(E)$ is a finite extension of $K'$. If both $E$ and $E'$ 
contain a primitive $p^M$-th root of unity then 
the groups $\Gamma _K^{\ab }/p^M$ and $\Gamma _{K'}^{\ab }/p^M$ are provided 
with the $P$-topological structure, cf. n.8.2, and the induced isomorphism 
$$f^*_M:\Gamma _{K'}^{\ab }/p^M\longrightarrow\Gamma _K^{\ab }/p^M$$ 
is $P$-continuous. 

Consider the inverse statement.

\proclaim{Theorem 6} With the above notation 
suppose that $p\geqslant 3$ and 
$\iota :\Gamma _{K'}\longrightarrow\Gamma _K$ is an isomorphism 
of profinite groups such that 
\newline 
{\rm a)} for all $j\in J(N)$, $\iota (\Gamma _{K'}^{(j)})=\Gamma _K^{(j)}$;
\newline 
{\rm b)} if $E,E'$ are finite extensions of $K$ and, resp., $K'$ such that 
the both contain $\zeta _{p^M}$, then the induced 
isomorphism 
$$\iota _M:\Gamma _{E'}^{\ab }/p^M\longrightarrow \Gamma _E^{\ab }/p^M$$
is $P$-continuous. 

Then there is a (unique) field isomorphism 
$f:\Bbb C_p(N)\longrightarrow\Bbb C_p(N)$ such that $f(K)=K'$ 
and $f=\iota ^*$.
\endproclaim 

\remark{Remark} Modulo some technical details  
and notation this statement has been announced  in [Ab5].  
\endremark 

\demo{Proof} 

9.3.1. Notice first, that $\iota $ induces for $1\leqslant r\leqslant N$, 
the group isomophisms 
$\iota (r):\Gamma _{K'(r)}\longrightarrow\Gamma _{K(r)}$. 
All these isomorphisms are also compatible with the corresponding 
ramification filtrations. 

In particular, $\iota (1)$ is a compatible with ramification 
filtration isomorphism of the absolute Galois groups of 
1-dimensional local fields $K(1)$ and $K'(1)$. Therefore, 
by the 1-dimensional case of a local analogue of the 
Grothendieck conjecture, cf. n.8.1, $\iota (1)$ is 
induced by a field isomorphism $f(1):\Bbb C_p\longrightarrow \Bbb C_p$ such 
that $f(1)(K(1))=K'(1)$. 
\medskip 

9.3.2. Prove the existence of $F_{\centerdot }, F_{\centerdot }'\in\Cal B^{fa}(N)$ 
such that for all $n\geqslant 0$, 
\newline 
a) $F_0\supset K$, $F_0'\supset K'$;
\newline 
b) $\iota (\Gamma _{F'_n})=\Gamma _{F_n}$;
\newline 
c) $\zeta _n\in F_n$ and $\zeta _n\in F'_n$, where $\zeta _n$ is 
a primitive $p^n$-th root of unity.
\medskip 

Let $E_0=\Bbb Q_p\{\{\tau _N\}\}\dots \{\{\tau _2\}\}$ be a basic $N$-dimensional 
local field. Then $K$ and $K'$ are its finite extensions with induced $F$-structures. 
Consider $E_{\centerdot }\in\Cal B(N)$ such that for all $n\geqslant 1$, 
$E_n=E_0(\zeta _n,\root {p^n}\of {\tau _2},\dots ,\root {p^n}\of {\tau _n})$. 
Clearly, $E_{\centerdot }\in\Cal B^a(N)$ 
(even more, $E_{\centerdot }\in\Cal B^{fa}(N)$). 

Let $L_{\centerdot }=KE_{\centerdot }$.  Then $L_{\centerdot }\in\Cal B^a(N)$ by Prop.4.1.  
Introduce $L_{\centerdot }'=\{L_n'\ |\ n\geqslant 0\}\in\Cal B(N)$ such that 
$\iota (\Gamma _{L'_n})=\Gamma _{L_n}$. Then $L_{\centerdot }'\in\Cal B^a(N)$ 
because $\iota $ is compatible with ramification filtrations. 

Suppose $n^*=n^*(L_{\centerdot })$ is the parameter for $L_{\centerdot }$ introduced in n.4.2.  
Clearly, $n^*$ can be taken  also as a  parameter for $L_{\centerdot }'$. 
Choose a finite extension 
$M(N-1)$ of   
$L_{n^*}(N-1)$ such that if $M=L_{n^*}M(N-1)$ then 
$(M, M(N-1))\in\LC (N)$ is standard, cf. Theorem 1. If necessary, we can enlarge 
$M(N-1)$ to satisfy the following property: if $M(N-1)'$ is such that 
$\iota (N-1)(\Gamma _{M(N-1)'})=\Gamma _{M(N-1)}$ and 
$M'=L'_{n^*}M(N-1)'$ 
then $(M',M(N-1)')\in\LC (N)$ is standard. Therefore, the towers 
$M_{\centerdot }=ML_{\centerdot }$ and $M_{\centerdot }'=L'_{n^*}M(N-1)'$ are such that 
for all $n\geqslant 0$, $\iota (\Gamma _{M_n'})=\Gamma _{M_n}$ and 
\linebreak 
$(M_{n^*},M_{n^*}(N-1)), (M'_{n^*},M'_{n^*}(N-1))\in\LC (N)$ are standard. 

Apply 
the above procedure to $(N-1)$-dimensional towers 
$M_{\centerdot }(N-1),M_{\centerdot }'(N-1)\in\Cal B^a(N-1)$ with  
a parameter $m^*\geqslant n^*$ and so on. Finally, we obtain 
finite separable extensions $F_{\centerdot }$ and $F'_{\centerdot }$ 
of $L_{\centerdot }$ and, resp., $L_{\centerdot }'$, which still satisfy 
the above requirements a)-c) but are already objects of the category 
$\Cal B^{fa}(N)$. 
\medskip 

9.3.3. Let $\Cal F\in X(F_{\centerdot })$ and 
$\Cal F'\in X(F_{\centerdot }')$, cf. section 5.  By Theorem 2 the group 
isomorphism  
$\iota $ induces the identification 
$$\iota _{F_{\centerdot }}:\Gamma _{\Cal F'}\longrightarrow \Gamma _{\Cal F}$$
which is compatible with ramification filtrations on these groups. 

Suppose finite extensions $\Cal E/\Cal F$ and $\Cal E'/\Cal F'$ 
are such that $\iota _{F_{\centerdot }}(\Gamma _{\Cal E'})=\Gamma _{\Cal E}$. 
If $\Cal E$ and $\Cal E'$ have standard $F$-structures then 
$\Cal E\in X(E_{\centerdot })$ and $\Cal E'\in X(E_{\centerdot }')$, 
where $E_{\centerdot }, E_{\centerdot }'\in\Cal B^{fa}(N)$ are 
finite separable extensions of $F_{\centerdot }$ and $F_{\centerdot }'$, 
respectively. Therefore, we can apply  Theorem 4 to deduce 
from the condition b) of the statement of our theorem  that 
for any $M\in\Bbb N$, 
the induced identification 
$$\iota _{F_{\centerdot }M}:\Gamma _{\Cal E'}^{\ab }
/p^M\longrightarrow \Gamma ^{\ab }_{\Cal E}/p^M$$
is $P$-continuous. 

Therefore, by the charactersitic $p$ case of the Grothendieck Conjecture,  
cf. Theorem 5 in n.9.2, the isomorphism $\iota _{F_{\centerdot }}$ is induced by 
a field isomorphism $f_R:R_0(N)\longrightarrow R_0(N)$ such that 
$f_R(\Cal R(\Cal F))=\Cal R(\Cal F')$.
\medskip 

9.3.4. Clearly, $f_R|_{R_0(1)}$ is induced by the $f(1):\Bbb C(N)_p\longrightarrow\Bbb C(N)_p$ 
from n.9.3.1. Therefore, $f_R$ leaves invariant the subgroup of Fontaine's elements 
$<\varepsilon>$ and by Lemma 7.4, 
$f_R$ is induced by a field automorphism $f:\Bbb C(N)_p\longrightarrow\Bbb C(N)_p$. 

The characteristic property of the field automorphism $f_R$ is that it 
transforms the action of any $\tau\in\Gamma _{K'}$ on $R_0(N)$ into the action of 
$\iota (\tau )\in\Gamma _K$ on $R_0(N)$. Therefore, 
$f$ satisfies the same property and we have 
$$f(K)=f(\Bbb C(N)_p^{\Gamma _{K}})=\Bbb C(N)_p^{\Gamma _{K'}}=K'.$$

So, $f\in\Hom _{\LF _0(N)}(K,K')$ and Theorem 6 is proved.
\enddemo

\ref\key
\by 
\paper 
\jour 
\vol 
\issue 
\yr 
\page 
\endref

\ref\key
\by 
\paper 
\jour 
\vol 
\issue 
\yr 
\page 
\endref

\ref\key
\by 
\paper 
\jour 
\vol 
\issue 
\yr 
\page 
\endref

\enddocument

\enddocument

\enddocument

\enddocument

\enddocument